\documentstyle[amssymb]{article}

\newcommand{\RR}{{\Bbb R}}
\newcommand{\p}{\frac{n+2}{n-2}}
\newcommand{\pp}{\frac{2n}{n-2}}
\newcommand{\ppp}{\frac{4}{n-2}}
\newcommand{\q}{\frac{2-n}{2}}
\newcommand{\e}{\varepsilon}
\newcommand{\al}{\alpha}
\newcommand{\Lap}{\Delta}
\newcommand{\del}{\partial}
\newcommand{\ML}{{\cal M}_\Lambda} 
\newcommand{\Mk}{{\cal M}_k}
\newcommand{\SL}{S^n \setminus \Lambda}
\newcommand{\La}{\Lambda}
\newcommand{\Li}{{\cal L}}
\newcommand{\dive}{{\rm div\,}}
\newcommand{\G}{{\cal G}}
\newcommand{\ha}{\frac{1}{2}}
\newcommand{\th}{\theta}
\newcommand{\be}{\bar{\e}}

\newcommand{\ba}{\bar{a}}

\begin{document}

\title{Refined asymptotics for constant scalar curvature metrics with 
isolated singularities}
\date{ } 
\author{Nick Korevaar\thanks{Supported by NSF grant \#DMS-9208666} \\ 
University of Utah \\ korevaar\@@math.utah.edu
\and Rafe Mazzeo\thanks{Supported by the NSF with a Young Investigator 
Fellowship and by grant \#DMS-9626382}\\ Stanford University \\ 
mazzeo\@@math.stanford.edu \and Frank Pacard \\ Universit\'e  Paris XII \\ 
Frank.Pacard\@@univ-paris12.fr \and Richard Schoen\thanks{Supported
by the NSF grant \#DMS-9504441 and the Guggenheim Foundation} \\ Stanford University 
\\ schoen\@@math.stanford.edu}  
\maketitle

\newtheorem{theorem}{Theorem}
\newtheorem{lemma}{Lemma}
\newtheorem{proposition}{Proposition}
\newtheorem{corollary}{Corollary}
\newtheorem{definition}{Definition}
\newtheorem{remark}{Remark}

\begin{abstract} We consider the asymptotic behaviour of positive solutions 
$u$ of the conformal scalar curvature equation, $\Delta u + \frac{n(n-2)}{
4}u^{\frac{n+2}{n-2}} = 0$, in the neighbourhood of isolated singularities in 
the standard Euclidean ball. Although asymptotic radial symmetry for such 
solutions was proved some time ago, \cite{CGS}, we present a much simpler 
and more geometric derivation of this fact. We also discuss a refinement, 
showing that any such solution is asymptotic to one of the deformed 
radial singular solutions. Finally we give some applications of these refined 
asymptotics, first to computing the global Poho\^zaev invariants of solutions 
on the sphere with isolated singularities, and then to the regularity of the 
moduli space of all such solutions.
\end{abstract}

\section{Introduction}
The problem we consider in this paper is to derive asymptotics for positive 
solutions of the conformally invariant semilinear elliptic equation 
\begin{equation}
\Lap u + \frac{n(n-2)}{4}u^{\p} = 0
\label{eq:1.1}
\end{equation}
which are defined in the punctured ball $B^n \setminus \{0\}$, and which are 
singular at the origin. It is well-known that a solution $u$ of this equation 
corresponds to a conformally flat metric 
\[
g = u^{\ppp}\delta
\]
which has constant scalar curvature $R_g = n(n-1)$.  Equation (\ref{eq:1.1}) 
is a special case of the more general equation relating the scalar curvatures 
of any two conformally related metrics. More specifically, if $g$ and $g' = 
u^{\ppp}g$ are any two such metrics, with corresponding scalar curvature 
functions $R_g$ and $R_{g'}$, respectively, then
\begin{equation}
\Lap_g u - \frac{n-2}{4(n-1)}R_gu + \frac{n-2}{4(n-1)}R_{g'}u^\p = 0.
\label{eq:1.3}
\end{equation}
The linear part of the operator on the left 
\begin{equation}
L_g = \Lap_g - \frac{n-2}{4(n-1)}R_g
\label{eq:1.31}
\end{equation}
is known as the {\it conformal Laplacian} associated with the metric $g$. 
It satisfies the conformal equivariance property that 
\begin{equation}
L_{g} (u\phi ) = u^{\frac{n+2}{n-2}}L_{g'}\phi
\label{eq:1.32}
\end{equation}
for any $\phi \in {\cal C}^\infty$ and any two metrics $g$ and $g'$ related 
as above. 

Asymptotics of solutions of equation (\ref{eq:1.1}) were proved by L.A. 
Caffarelli, B. Gidas and J. Spruck \cite{CGS} using a rather complicated 
version of Alexandrov reflection. Their result states that an arbitrary 
positive solution $u$ of (\ref{eq:1.1}) satisfies
\begin{equation}
u(x) = |x|^{\q}v_\e(-\log |x| + T)(1+o(1)),
\label{eq:1.4}
\end{equation}
where $u_\e(x) = |x|^{\q}v_\e(-\log |x|)$ is one of a one-parameter family 
of radial solutions of (\ref{eq:1.1}). We call the $v_\e$ and the $u_\e$ 
{\it Fowler solutions} and discuss them at greater length in the next 
section ({\it Fowler solutions} are also called {\it Delaunay type 
solutions} in many works which we will refer to).  In recent work 
\cite{L}, C.C. Chen and C.S. Lin have generalized this result and slightly 
simplified their argument, although their proof is still quite technical. 
This asymptotics theorem is related to a simpler result due to P. Aviles 
\cite{A} concerning the subcritical equation, where the exponent $\p$ is 
replaced by a smaller number $p$ in the range $[\frac{n}{n-2}, \p)$. 
P. Aviles shows 
that an arbitrary positive solution of the subcritical equation with an 
isolated singularity at the origin converges, in the sense of (\ref{eq:1.4}),
to a radial solution of the same equation.

The present work is a combination of two separate projects which are 
sufficiently closely related that it has seemed natural to publish them 
together. The first, undertaken several years ago by the first and fourth 
authors (K \& S), is to supply an alternate argument, more geometric and 
much simpler than that of \cite{CGS}, to establish (\ref{eq:1.4}). This 
also gives a slightly stronger estimate, improving the $o(1)$ remainder term 
to $O(|x|^{\al})$ for some $\al > 0$. The second project, quite recently 
undertaken by the second and third authors (M \& P) and inspired by their 
recent construction \cite{MP2}, concerns a refinement of these asymptotics, 
in some sense obtaining a second term in the expansion, along with various 
applications for these refined asymptotics. More precisely, it is possible 
to improve (\ref{eq:1.4}) using, instead of $v_\e$, a family of deformations 
of the Fowler solutions, parametrized by vectors $a \in \RR^n$, 
which we denote $v_{\e,a}$ and which are defined by 
\[
v_{\e,a}(t,\theta) = |\theta - ae^{-t}|^\q v_\e( t + \log|\theta - ae^{-t}|),
\]
where $\theta = x/|x|$. For all $a \in \RR^n$ and $T \in \RR$ we define
\[
u_{\e,a, T}(x) = |x|^{\q} v_{\e,a}( -\log|x|+ T,x/|x|).
\]
The main result is

\begin{theorem} Let $u$ be an arbitrary weak solution of the equation
(\ref{eq:1.1}) which is positive on the ball $B^n(0) \subset \RR^n$ and in 
${\cal C}^\infty(B^n(0)\setminus\{0\})$. Then either $u \in 
{\cal C}^\infty(B^n(0))$, or else there exists some choice of parameters 
$\e \in (0, \e_0]$ ($\e_0 = ((n-2)/n)^{(n-2)/4}$), $a \in \RR^n$ and some 
constants $T \in \RR$ and $\al > 1$ such that
\[
u(x) = |x|^{\q}(v_{\e,a}(-\log|x| + T,x/|x|) + O(|x|^\al)), \quad \mbox{ as }
|x| \rightarrow 0.
\]
Alternatively,
\[
u(x) = u_{\e,a,T}(x) + O(|x|^\beta) \quad \mbox{ as } |x| \rightarrow 0
\]
for some constant $\beta > \frac{4-n}{2}$. In particular, there are upper and 
lower bounds for $u$ of the form $C_1 |x|^\q \le u \le C_2 |x|^\q$, with 
$0 < C_1 \le C_2$, and the corresponding metric $g = u^\p \delta$ is complete 
at $0$.
\label{th:refas}
\end{theorem}

Solutions of equation (\ref{eq:1.3}) on the sphere $S^n$, singular along a 
closed set $\Lambda$, giving conformal factors relating the standard metric 
to a new metric of scalar curvature $n(n-1)$ complete on $S^n \setminus 
\Lambda$, were first considered in the work \cite{SY} of the fourth author 
and S.T. Yau. For solutions to exist it is necessary that the Hausdorff 
dimension of $\Lambda$ be less than or equal to $(n-2)/2$. Under mild 
geometric hypotheses, these solutions are shown there to extend to global 
weak solutions on $S^n$. Around the same time, solutions with isolated 
singularities (and with some more complicated singular sets arising 
essentially as limit sets of Kleinian groups) on $S^n$ were constructed 
by one of us \cite{S1}. Since then, solutions singular along an 
arbitrary disjoint 
collection of submanifolds, and giving metrics of {\it positive} scalar 
curvature, have been constructed in a succession of works, \cite{MS}, 
\cite{P}, \cite{MP1}, \cite{MP2}. (We do not discuss the by now extensive 
set of results on the analogous problem for metrics of negative scalar 
curvature.) 

Emerging from \cite{MP1} and \cite{MP2} (as well as \cite{S1}) is the lesson 
that solutions with isolated singularities are far more rigid objects than 
solutions with higher dimensional singular sets. This is apparent if one 
fixes the singular set $\Lambda$ and studies the moduli space $\ML$ of all 
solutions singular along $\Lambda$. We shall henceforth restrict attention 
to the case where $\Lambda$ is a finite disjoint collection of submanifolds,
of (possibly varying) dimensions less than or equal to $(n-2)/2$. When at 
least one component of $\Lambda$ is of positive dimension, $\ML$ is infinite 
dimensional and is modelled locally on a Banach manifold. Closely related 
to Theorem~\ref{th:refas} is a result on refined asymptotics, analogous to the 
one here but for solutions singular along higher dimensional submanifolds, 
given in \cite{M}. One sees from these asymptotics that solution metrics 
need not have bounded curvature, and that the correct way to parametrize the
moduli space $\ML$ then is via an asymptotic Dirichlet problem. On the other 
hand, when $\Lambda$ is a set of $k$ points on the sphere, then $\ML$ is 
finite dimensional.  It is proved in \cite{MPU1} that in this case, $\ML$ 
is even a locally real analytic set with virtual dimension $k$.
In \cite{MP2} it is shown that for generic configurations $\Lambda$, $\ML$ 
attains this dimension. The sharper asymptotics proved here have applications 
to the study of $\ML$ when $\Lambda$ is finite. In particular, corresponding 
to any conformal Killing vector field on the sphere one can define a 
Poho\^zaev invariant of the (global) solution, depending on the parameters 
$\e$ and $a$ associated to each singular point. Using the refined asymptotics 
one can evaluate these invariants explicitly; one obtains a finite number 
of real analytic equations satisfied by the set of parameters associated to 
each solution. While these equations are not enough to determine the moduli 
space completely, they do lead to new information about it. This is done 
in \S 6. There are two other applications of the refined asymptotics 
considered here. Then, in \S 7, we discuss the nondegeneracy of the 
moduli space $\ML$ near its ends. 

For any $k > 0$, Theorem~\ref{th:refas} gives the existence of a 
`parameter map' $\Phi$ on the space of solutions defined
on all of ${\RR}^n$, with exactly $k$ singularities; $\Phi$ 
associates to a solution its set of singular points $p_1, \ldots, 
p_k$ and the associated parameters $T_j, a_j$ and $\e_j$ at each 
$p_j$. $\Phi$ is defined more carefully in (\ref{eq:10.1}) below.
This map should be quite useful in the global study of the moduli 
spaces $\ML$ when $\Lambda$ is finite. It seems quite likely
that $\Phi$ is injective, but the proof of this fact may be
subtle. We leave this as an interesting open problem. 

We conclude this introduction by remarking on the relationship of the results 
here with what is known for noncompact constant (nonzero) mean curvature 
surfaces $\Sigma \subset \RR^3$ which are (Alexandrov) embedded and have $k$ 
ends. Although seemingly quite different, this extrinsic geometric problem 
bears many similarities to the intrinsic conformal scalar curvature equation. 
Indeed, all the results known in one setting seem to have direct analogues 
in the other setting, although the proofs are sometimes different. Any end 
of such a surface $\Sigma$ is asymptotic to one of the Delaunay surfaces 
(these are a family of rotationally symmetric CMC surfaces in $\RR^3$); 
this follows from work of the 
first author, Kusner and Solomon, \cite{KKS}, and Meeks, \cite{Me}. The 
deformations of the radial solutions used here, and the corresponding refined 
asymptotics, corresponds in the CMC setting to translations and rotations of 
the asymptotic axes. More analytically, the difference is that the Jacobi 
fields corresponding to these deformations, in the CMC context, are at most 
linearly growing, while in the scalar curvature setting they are exponentially
growing or decaying. Amongst the various applications of the refined 
asymptotics we develop here, the global balancing formula of \S 6 corresponds 
to the computationally more simple balancing formula in the CMC context 
discussed, for example, in \cite{KKS}. The nondegeneracy argument of 
\S 7 has a direct analogue in the CMC setting, and to our knowledge, is 
new there as well. 

\section{Preliminaries} 
In this section we define and discuss in some detail the $(n+1)$-parameter 
family of functions $v_{\e,a}$ which arise in the statement of 
Theorem~\ref{th:refas}. We also study the linearization of (\ref{eq:1.1}) 
about any one of the the Fowler solutions $v_{\e}$, since this will 
be a crucial tool in the later analysis. Parts of this material can also be 
found in \cite{S1} and \cite{MPU1}. 

In the ensuing discussion it will be convenient to use the conformal 
equivariance of the problem to express (\ref{eq:1.1}) relative to different 
background metrics in the flat conformal class. The three most natural 
choices are the standard Euclidean metric, the cylindrical metric $dt^2 + 
d\theta^2$ on $\RR \times S^{n-1}$, and the standard spherical 
metric on $S^n$. Equation (\ref{eq:1.1}) is already expressed relative to 
the Euclidean metric. Relative to the cylindrical metric this equation becomes
\begin{equation}
\frac{d^2v}{dt^2} + \Lap_{\theta}v - \frac{(n-2)^2}{4}v + \frac{n(n-2)}{4}
v^\p = 0.
\label{eq:2.1}
\end{equation}
The transformation between (\ref{eq:1.1}) and (\ref{eq:2.1}) can also be seen 
directly by setting 
\[
u (x)= |x|^\q v(-\log |x|, x/|x| )
\]
and choosing a new independent variable $t = -\log |x|$. We shall not need 
the explicit form of this equation relative to the spherical metric. 

At each stage below we shall use whichever of these ambient metrics is most 
convenient for the computations at hand.

\subsection{Fowler solutions}
We now define the special family of solutions mentioned in the introduction 
and develop their relevant properties. We start by considering the radial
solutions of (\ref{eq:1.1}) on $\RR^n$ which are singular at the origin. For 
this problem it is simplest to use the cylindrical metric and equation 
(\ref{eq:2.1}). 

The first result characterizes the global solutions of (\ref{eq:2.1}).
\begin{proposition}
Let $v$ be any positive solution of (\ref{eq:2.1}) defined on the whole 
cylinder $\RR \times S^{n-1}$. Then $v$ is independent of the spherical 
variable $\theta$, hence depends only on the variable $t$. Moreover we always 
have $ 0 < v(t) \leq 1$ for all $t \in \RR$. 
\label{pr:2.1}
\end{proposition}
The proof of this uses an Alexandrov reflection argument, and is contained in 
\cite{CGS} and also in \cite{S2}. 

Hence if $v$ is a global solution, then it satisfies the autonomous ODE 
obtained by dropping the $\theta$ differentiation in (\ref{eq:2.1}). We 
analyze this by converting it to a system of first order equations on the 
phase plane. Namely, setting $w = \del_t v$, the equation becomes
\[
\frac{d\,}{dt} v = w, \qquad \frac{d\,}{dt} w = \frac{(n-2)^2}{4} v - 
\frac{n(n-2)}{4}v^{\frac{n+2}{n-2}}.
\]
This system is Hamiltonian, with corresponding energy function
\begin{equation}
H(v,w) = w^2 - \frac{(n-2)^2}{4}v^2 + \frac{(n-2)^2}{4}v^\pp.
\label{eq:2.2}
\end{equation}
In particular, if $v(t)$ is a solution of (\ref{eq:2.1}) then the path 
parametrized by $(v(t),v'(t))$ is contained within a level set of $H$, hence 
$H(v,v')$ does not depend on $t$. There is a homoclinic trajectory lying in 
the level set $\{H = 0\}\cap \{v > 0\}$, and a one-parameter family of closed 
level sets contained in the bounded set ${\cal O} \equiv \{H < 0\} \cap 
\{v > 0\}$. All other trajectories and level sets pass into the region where 
$v \le 0$, so we do not consider them here. The compact level sets in 
${\cal O}$ correspond to a family of periodic solutions of (\ref{eq:2.1}) 
which we shall denote $v_\e(t)$. The parameter $\e$ denotes the minimum value 
attained by the solution, and we will often refer to it as the {\it necksize} 
of the solution. It varies in the interval $(0,\e_0]$, where 
$\e_0^{\frac{4}{n-2}} = \frac{n-2}{n}$. The level set of $H$ corresponding 
to $v_\e$ is
\begin{equation}
H \equiv H(\e) = \frac{(n-2)^2}{4} \big(\e^\pp - \e^2 \big).
\label{eq:hamenge}
\end{equation}

Before continuing with the rest of our treatment of this special family of
solutions we digress slightly to discuss the energy of arbitrary (nonnegative) 
solutions of (\ref{eq:2.1}). Thus let $v$ be any such solution. Multiply
this equation by $\del_t v$; a small calculation yields
\[
\frac{1}{2}\del_t \big(\del_t v \big)^2 + \dive_\theta\big( \del_t v 
\nabla_\theta v \big) -\frac{1}{2} \del_t |\nabla_\theta v |^2 - 
\frac{(n-2)^2}{8} \del_t (v^2) + \frac{(n-2)^2}{8}\del_t (v^\pp) = 0.
\]
We can rephrase this by noting that it is equivalent to the vanishing of the 
divergence of the vector field 
\[
W = \left(\frac{1}{2} (\del_t v)^2 -\frac{1}{2}  |\nabla_\theta v |^2 
- \frac{(n-2)^2}{8} v^2 + \frac{(n-2)^2}{8} v^\pp \right)\del_t +  (\del_t v) 
\nabla_\theta v \nabla_\theta.
\]
Hence the integral of $\langle W, \del_t\rangle$ on the sphere $t = T$ (with
respect to the volume form for the metric $g$) gives a number which, by the 
divergence theorem, does not depend on $T$. Thus we obtain an invariant of 
the solution $v$, which we shall call the `radial Poho\^{z}aev invariant' of 
$v$ and denote ${\cal P}_{\rm ra}(v)$. It is a special case of the more 
general Poho\^{z}aev invariants we shall discuss later. Note that the radial 
Poho\^{z}aev invariant of the Fowler solution $v_\e$ is simply 
$\frac{1}{2}\omega_{n-1}H(\e)$, where $\omega_{n-1}$  is the volume of 
$S^{n-1}$ and $H(\e)$ is the Hamiltonian energy (\ref{eq:2.2}). 

We now return to the discussion of the special family of solutions. There are 
actually two parameters for these periodic solutions. The first is this 
necksize parameter $\e$, while the second corresponds to the value of the 
solution at $t = 0$ and we will denote by $T_\e$ the period of $v_\e$. We 
normalize the functions $v_\e$ by assuming that $v_\e(0) = \min v_\e$. Hence 
the complete family of positive radial solutions is $v_\e(t + T)$ for 
$T\in {\RR}$ and $\e \in (0, \e_0]$.  The two extremes of this family are when
$v$ is the constant solution $v \equiv \e_0$ (geometrically, this corresponds 
to the correct magnification of the cylinder which has scalar curvature 
$n(n-1)$), and the limits $v_\e(t) \longrightarrow 0$ and $v_\e(t+ \ha T_\e) 
\longrightarrow (\cosh t)^{\q}$ as $\e \rightarrow 0$. This last explicit 
solution corresponds to the conformal factor $(\cosh t)^{-2}$ transforming the 
cylinder to the punctured sphere $S^n \setminus \{p,-p\}$. 

It is known that $T_{\e}$, the period of $v_\e$, is monotone in $\e$, 
converging to $2 \pi / \sqrt{n-2}$ as $\e \rightarrow \e_0$ and increasing to 
$\infty$ as $\e \rightarrow 0$. Geometrically, however, the metrics 
$g_\e \equiv v_\e^\p (dt^2 + d\theta^2)$ converge to a bead of spheres 
arranged along a fixed axis.  This family of radial solutions interpolates 
between the cylinder and this singular limit. 

We may also transform these back to solutions of (\ref{eq:1.1}), to obtain the
 solutions 
\[
e^{\q T} u_\e(e^{-T}x) = |x|^\q v_\e(-\log|x| + T)
\]
on $\RR^n \setminus \{0\}$. 

As indicated earlier, we shall consider, in addition, an $n$ parameter 
family of deformations of these radial solutions. Geometrically these may
be understood as follows. Using the Euclidean background metric, 
and starting with a radial solution $u$, we first take the Kelvin transform 
$\tilde{u} (x) = |x|^{2-n}u(x/|x|^2)$. It is not difficult to check
that $\tilde{u}$ is also a solution of (\ref{eq:1.1}). Then translate
$\tilde{u}$ by some vector $a \in \RR^n$, and finally take the
Kelvin transform once again to obtain $u_{\e,a}$. Note that this is
the pullback by the composition of three conformal transformations,
first inversion in the unit sphere $\{|x| = 1\}$, then Euclidean translation,
which is a parabolic transformation fixing infinity, and then inversion
once again. Each of these elementary transformations takes the space of 
solutions of (\ref{eq:1.1}) to itself. Using instead the spherical metric 
as the background, if $u$ is transformed to a solution on $S^n$ 
which is singular
at two antipodal points $\{p,-p\}$, then these steps correspond
to first reflecting across the equator determined by $p$ and $-p$, 
then applying a parabolic transformation $F$ which fixes the point $p$ and
carries another point $q$ to $p$, and finally reflecting across the
equator again. Finally, note that relative to the cylindrical background
metric, the final and initial transformations correspond even more 
simply to the reflection $t \rightarrow -t$. The parabolic
translation is more complicated to describe on the cylinder.

Now let us perform these steps explicitly 
on $u_{\e}(x) = |x|^\q v_{\e}(-\log |x|)$. 
After the first inversion, we obtain $|x|^\q v_{\e}( \log |x|)$.
Translation carries this to $|x-a|^\q v_{\e} (\log |x-a|)$. Performing
the final inversion yields at last
\[
u_{\e,a}(x) = |\frac{x}{|x|^2} - a|^\q v_\e (\log|\frac{x}{|x|^2} - a|)
= |x|^\q |\theta - a|x||^\q v_\e( -\log|x| + \log|\theta - a|x||),
\]
where $\theta = x/|x|$. We can also define the associated functions
\[
v_{\e,a}(t,\theta) = |\theta - ae^{-t}|^\q v_\e( t + \log|\theta - ae^{-t}|)
\]
on the cylinder; these are solutions of (\ref{eq:2.1}).  Notice that 
these functions are regular except at $t = \log |a|, \theta = a/|a|$,
and in particular are smooth for $t$ sufficiently large. 

We conclude this section by calculating the expansion of $u_{\e,a}$
near $|x| = 0$. To do this we first observe that
\[
\big|\frac{x}{|x|} - a|x|\big|^{\q} = 1 + \frac{n-2}{2} a \cdot x + O(|x|^2).
\]
Similarly we calculate
\[
\log \big|\frac{x}{|x|} -a|x| \big| = - a \cdot x + O(|x|^2),
\]
and so
\[
v_{\e}(-\log |x| - a \cdot x + O(|x|^2)) = v_\e(-\log |x|) -
v_{\e}'(-\log |x|) (a \cdot x) + O(|x|^2).
\]
Putting these altogether we finally obtain
\begin{eqnarray}
u_{\e,a}(x) = |x|^\q\left(v_\e(-\log |x|) + a\cdot x 
\big( -v_\e'(-\log |x|)+\frac{n-2}{2}v_\e(-\log |x|)\big) + O(|x|^2) \right) 
\nonumber \\
= u_\e(x) + |x|^\q (a \cdot x) (-v_\e' + \frac{n-2}{2}v_\e) + O(|x|^{
\frac{6-n}{2}}).
\label{eq:2.25}
\end{eqnarray}
In particular, $|v_{\e,a}(t,\theta) - v_\e(t)| \le Ce^{-t}$. 

We have now obtained the families of solutions $u_{\e,a}(x)$
and $v_{\e,a}(t,\theta)$ for the problems (\ref{eq:1.1}) and
(\ref{eq:2.1}), respectively. 
For all $T\in {\Bbb R}$, changing $v_{\e}(t)$ into $v_{\e}(t+T)$ in the
construction above leads to families of solutions $u_{\e,a,T}(x)$
and $v_{\e,a,T}(t,\theta)$ for the problems (\ref{eq:1.1}) and
(\ref{eq:2.1}), respectively.
We note again that the parameters
$a \in \RR^n$ and $T \in \RR$ correspond to
explicit geometric (conformal) motions; only the parameter $\e$ does
not arise from an extrinsic motion. For simplicity, when $a = 0$, we shall 
denote
these functions simply by $u_\e$ and $v_\e$. 

\subsection{The linearized equation}
We shall now consider the linearization 
\begin{equation}
\Li_\e = \frac{\del^2\,}{\del t^2} + \Lap_{\theta} - \frac{(n-2)^2}{4} + 
\frac{n(n+2)}{4}v_\e^\ppp
\label{eq:2.3}
\end{equation}
of the operator in (\ref{eq:2.1}) around one of the {\it radial} Fowler 
solutions $v_\e$. Our primary interest is in the mapping properties of $\Li_\e$,
which we shall review from \cite{MPU1} and \cite{MP2}. 

This operator has periodic coefficients, hence may be studied by classical 
Flocquet - or Bloch wave - theoretic methods, as in \cite{MPU1}, but also by 
separation of variables and elementary ODE methods, as in \cite{MP2}. 
We shall need to refer to results derived from both of these methods. 

Let $\{\lambda_j,\; \chi_j(\theta)\}$ be the eigendata of $\Lap_{S^{n-1}}$. 
We use the convention that the eigenvalues are listed with  multiplicity,
so that $\lambda_0 = 0$, $\lambda_1 = \ldots = \lambda_{n}= n-1$, 
$\lambda_{n+1} = 2n$, etc.
Then $\Li_\e$ decouples into infinitely many ordinary differential operators
\[
\Li_{\e,j} = \frac{d\,^2}{dt^2} + \big( \frac{n(n+2)}{4}v_\e^\ppp -
\frac{(n-2)^2}{4} - \lambda_j \big).
\]
When $j > n+1$, the term of order zero in each of the $\Li_{\e,j}$ is negative
because $\lambda_j \geq 2n$ and $v_\e < 1$; hence for these values
of $j$, $\Li_{\e,j}$ satisfies the maximum principle. The same conclusion
holds for $\Li_{\e,j}$, $j = 1, \ldots, n$, because conjugating 
by an appropriate power of $v_{\e}$ 
yields an operator also with negative term of order zero.
It follows from these facts that the 
$L^2$ spectrum of $-\Li_{\e,j}$ is contained in $(0,\infty)$ when $j>0$. 
On the other hand, $0$ is contained in the {\it essential} 
spectrum of $-\Li_{\e,0}$. 
In particular, $\Li$ does not have closed range on $L^2(\RR \times S^{n-1};
dt\, d\theta)$, but the difficulty is localized to the ground eigenspace
of the Laplacian on the cross-section.

The failure of $\Li_\e$ to have closed range is caused by its Jacobi
fields, i.e. by the solutions of $\Li_\e \psi = 0$. To study these,
it suffices to consider the solutions of the induced problems
$\Li_{\e,j} \psi_j = 0$. It can be proved, cf. \cite{MPU1}, 
that for each $j > 0$ there
are two normalized linearly independent solutions $\psi_{\e,j}^\pm$ and 
an associated constant $\gamma_{\e,j}$ such that $|\psi_{\e,j}^\pm (t)| 
\le e^{\mp \gamma_{\e,j} t}$ for all $t$.
The solutions  $\psi_{\e, 0}^+(t)$, $\psi_{\e,0}^-(t)$ for $\Li_{\e,0}$ 
are bounded 
and linearly growing, respectively, and so it is natural to 
define $\gamma_{\e,0} = 0$. We immediately deduce the following 
\begin{proposition}
Suppose that $\Li_\e \psi = 0$ on $\RR^+ \times S^{n-1}$ and $\psi = 
O(e^{-\gamma t})$ as $t \rightarrow \infty$, for some $\gamma >0$. 
Then $\psi(t,\theta) = \sum_{j = j_0}^\infty \psi_{\e,j}^+(t)\chi_j(\theta)$,
where $j_0$ is the first integer such that
$\gamma_{j_0} > \gamma $.
\label{pr:2.2}
\end{proposition}

Somewhat remarkably, we may deduce the explicit forms of these solutions
when $j = 0, \ldots, n$. This rests on the fact that whenever we
have a one-parameter family of solutions of the nonlinear equation
(\ref{eq:2.1}), then the differential of this family gives
a solution of the linearized equation. In the previous subsection
we found several different one-parameter families of solutions.
The simplest are the families
\[
T \longrightarrow v_\e(t + T), \qquad \mbox{ and }\qquad \e \longrightarrow v_\e(t).
\]
Differentiation with respect to either of these two parameters
yields the solutions corresponding to the cross-sectional
eigenvalue $\lambda_0 = 0$. 
\begin{equation}
\psi_{\e,0}^+(t) \equiv \frac{d\,}{dT}v_\e(t+T) = v_\e'(t), \qquad
\psi_{\e,0}^-(t) \equiv \frac{d\,}{d\e}v_\e(t) = \dot{v}_\e(t).
\label{eq:2.35}
\end{equation}
We may also differentiate $v_{\e,a}$ with respect to $a_j$, using the 
expansion (\ref{eq:2.25}), to get 
\begin{equation}
\psi_{\e,j}^+(t)\chi_j(\theta) \equiv \frac{d\,}{da_j} v_{\e,a}(t, \theta) = 
e^{-t}\big(-v_\e'(t) + \frac{n-2}{2}v_\e(t)\big) \chi_j (\theta);
\label{eq:2.36}
\end{equation}
these give the exponentially decreasing solutions of of $\Li_{\e,j} \psi = 0$, 
$j = 1, \ldots, n$. Since these functions are all equal to one
another, it is often simpler to denote any one of them by just
$\psi_{\e,1}^+$.  In any case, using them we may restate (\ref{eq:2.25}) as
\begin{eqnarray}
u_{\e,a}(x) = |x|^\q\left( v_\e(-\log |x|) + |x|\big(\sum_{j=1}^{n} a_j
\chi_j(\theta) \psi_{\e,j}^+(-\log |x|)\big) + O(|x|^2) \right) \nonumber \\
= |x|^\q\left( v_\e(-\log |x|) + (a \cdot x) \psi_{\e,1}^+(-\log |x|) +
O(|x|^2)\right) .
\label{eq:2.37}
\end{eqnarray}
Notice that $\psi_{\e,1}^+= \ldots = \psi_{\e,n}^+$.
Although they are not required later, we also note that the exponentially 
increasing solutions $\psi_{\e,j}^-(t)$ of $\Li_{\e,j}$, $j = 1, \ldots, n$, 
are obtained from $\psi_{\e,j}^+(t)$ simply by replacing $t$ by $-t$. 

It follows from (\ref{eq:2.36}) that the `indicial root' $\gamma_{\e,j} 
\equiv 1$, $j = 1, \ldots, n$. It is somewhat unexpected that for
$j \le n$, $\gamma_{\e,j}$ is independent of $\e$; one can show however 
that all other $\gamma_{\e,j}$ depend nontrivially on $\e$, and in fact tend 
to $\displaystyle{\big( \frac{(n-2)^2}{4} +\lambda_j \big)^{1/2}}$
as $\e \rightarrow 0$ \cite{MP2}. 

The numbers $\{\pm \gamma_{\e,j}\}$ are analogous to the indicial roots
of a Fuchsian operator, and they influence the mapping properties
in exactly the same way. To state these mapping properties 
most suitably for our purposes we define weighted H\"older spaces on the 
(half-) cylinder
\[
\forall t_0\in {\RR}, \qquad 
{\cal C}^{k,\al}_{\gamma}( [t_0, +\infty) \times S^{n-1}) = \{ w = 
e^{t\gamma}\bar{w}:
\bar{w} \in {\cal C}^{k,\al}([t_0, +\infty) \times S^{n-1}) \}.
\]

If ${\cal F}$ is any space of functions defined on the half-cylinder 
$[t_0, +\infty) \times S^{n-1}$, 
we let ${\cal F}_{{\cal D}}$ denote the subspace of functions of ${\cal F}$
vanishing at $t = t_0$ for some fixed $t_0$. We cannot 
choose $t_0$ arbitrarily because it is necessary for the following result that
$\psi_j^+(t_0) \neq 0$ for $j \le n$; this condition is always fulfilled for 
$j=1, \ldots , n$ but fails for $j=0$ when $t_0 \in \{(\ell/2)T_\e: \ell 
\in {\Bbb Z}\}$), $T_\e$ being the period of $v_\e$.  We also define   
\[
E_{\e,0} = \mbox{Span} \{\psi_{\e,0}^{+}(t) \} \subset 
{\cal C}^{2,\al}_0 ( [t_0 , +\infty) \times S^{n-1})
\]
and
\[
E_{\e,1} = \mbox{Span} \{\psi_{\e,j}^+(t)\chi_j(\theta) : j= 1, \ldots , n\}
\subset {\cal C}^{2,\al}_{-1}( [t_0 , +\infty) \times S^{n-1}).
\]

Using this notation, the result we need is
\begin{proposition} Assume that $t_0 \neq (\ell/2)T_\e$, $\ell 
\in {\Bbb Z}$. Recalling that $\gamma_{\e,0} = 0$, $\gamma_{\e,1} = 
\ldots = \gamma_{\e,n} = 1$ and $\gamma_{\e,n+1} > 1$, if 
$0 < \gamma < 1$, then
\[
\Li_\e : \big[{\cal C}^{2,\al}_{-\gamma}( [t_0 , \infty) \times S^{n-1}) 
\oplus E_{\e,0} \big]_{\cal D} 
\longrightarrow {\cal C}^{0,\al}_{-\gamma}( [t_0 , \infty) \times S^{n-1}) 
\]
is a surjective Fredholm mapping. In particular, there exists a 
bounded right inverse
\[
\G_{\e,0} : {\cal C}^{0,\al}_{-\gamma}( [t_0 , \infty) \times S^{n-1}) 
\longrightarrow \big[ {\cal C}^{2,\al}_{-\gamma}( [t_0 , \infty) 
\times S^{n-1}) \oplus E_{\e,0} \big]_{\cal D}
\]
so that $\Li_\e \G_{\e,0} = I$. If $1 < \gamma < \gamma_{\e, n+1}$, then 
\[
\Li_\e : \big[{\cal C}^{2,\al}_{-\gamma}( [t_0 , \infty) \times S^{n-1}) 
\oplus E_{\e,0} \oplus E_{\e,1} \big]_{\cal D} 
\longrightarrow {\cal C}^{0,\al}_{-\gamma}( [t_0 , \infty) \times S^{n-1})
\]
is also a surjective Fredholm mapping, and there exists a bounded right
inverse
\[
\G_{\e,1}: {\cal C}^{0,\al}_{-\gamma}( [t_0 , \infty) \times S^{n-1}) 
\longrightarrow \big[ {\cal C}^{2,\al}_{-\gamma} ( [t_0 , \infty) 
\times S^{n-1})\oplus E_{\e,0} \oplus E_{\e,1} \big]_{\cal D}
\]
with $\Li_\e \G_{\e,1} = I$. 
\label{pr:map}
\end{proposition}

One may find solutions of $\Li_\e w = f$ for $f \in {\cal C}^{0,\al}_{
-\gamma}( [t_0 , +\infty) \times S^{n-1})$ by considering the induced 
differential operators $\Li_{\e,j}$ 
on each eigenspace of $\Delta_\theta$ separately. When $j > n$, 
solutions are easily obtained using the fact that $\Li_{\e,j}$ satisfies the 
maximum principle. This is done explicitly in \cite{MP2}. To obtain
solutions when $j \le n$ it is easiest to use the Fourier-Laplace
methods of \cite{MPU1}. Although these methods are $L^2$-based, it is
not hard to check that $w\in {\cal C}^{2,\al}_{-\gamma}( [t_0 , +\infty) 
\times S^{n-1})$ if
$f \in {\cal C}^{0,\al}_{-\gamma}( [t_0 , +\infty) \times S^{n-1})$. 
The necessity of adding $E_{\e,0}$ and $E_{\e,1}$ to the domains to obtain 
surjectivity for the homogeneous Dirichlet problem is a simple form
of the linear regularity theorem of \cite{MPU1}, but also follows
from easy ODE arguments. 

It is important in this result that the weight factor $-\gamma$
does not equal one of the $\pm \gamma_{\e,j}$. In fact, on any one of the
spaces ${\cal C}^{2,\al}_{\pm \gamma_{\e,j}}( [t_0 , +\infty) \times S^{n-1})$ 
the operator $\Li_\e$ does not have closed range. 

This proposition may be used not only to find solutions of
the equation $\Li_\e u = f$, but also to obtain decay properties
of solutions which are already given.
\begin{corollary} Suppose that $\mu$ and $\gamma$ are weight parameters
and $\Li_\e u = f$ where $u \in {\cal C}^{2,\al}_{-\mu}( {\Bbb R}^+
\times S^{n-1})$ and $f \in {\cal C}^{0,\al}_{-\gamma}( {\Bbb R}^+ 
\times S^{n-1})$. If $0 < \mu < \gamma < 1$, then $u \in 
{\cal C}^{2,\al}_{-\gamma}( {\Bbb R}^+ \times S^{n-1})$. If $0 < 
\mu < 1 < \gamma < \gamma_{\e,n+1}$, then $u \in {\cal C}^{2,\al}_{-\gamma}
( {\Bbb R}^+ \times S^{n-1}) \oplus E_{\e,1}$.
\label{cor:2.rinv}
\end{corollary}

The proof of this corollary is straightforward. For example,
if $0 < \mu < \gamma < 1$, fix some $t_0>0$ satisfying the assumptions of 
Proposition~\ref{pr:map}, then $v +c \psi^+_{\e,0} = \G_{\e,0} f \in 
[{\cal C}^{2,\al}_{-\gamma} ( [t_0 , +\infty) \times S^{n-1}) \oplus 
E_{\e,0}]_{\cal D}$ is also a solution of $\Li_\e v = f$. 
Hence $h = u - v$ is a homogeneous solution, $\Li_\e h = 0$, which
is also exponentially decaying as $t \rightarrow +\infty$. This means 
that $h \in {\cal C}^{2,\al}_{-1} ( [t_0 , +\infty) \times S^{n-1})$, 
and so we conclude that $u$ has the stated decay. The other cases 
are proved similarly. 

\section{Upper and lower bounds for solutions}

In this section we derive upper and lower bounds on solutions defined in the
punctured unit ball. We further show that any such solution has a nonzero
radial Poho\^zaev invariant ${\cal P}_{\rm ra}$. Specifically we prove the 
following.
\begin{theorem}
Assume that $u$ is a nonnegative smooth solution of (\ref{eq:1.1}) 
defined in the punctured unit ball $B^n\backslash \{0\}$. Either $u$ 
extends as a smooth solution to
the ball, or there exist positive constants $C_1$, $C_2$ such that
$$ 
C_1|x|^{(2-n)/2}\leq u(x)\leq C_2|x|^{(2-n)/2}.
$$
Furthermore, the radial Poho\^zaev invariant of $u$ is nonzero.
\label{th:3.1}
\end{theorem}
As stated, the constants $C_1$ and $C_2$ depend on the solution $u$. We will
see in this section that the upper bound $C_2$ follows from a more precise 
universal upper bound, while the discussion of the previous section, and
in particular the fact that the infimum of $v_\e$ tends to zero with
$\e$, shows that the lower bound depends in a more delicate way on the 
solution. It seems
to be unknown how one should express the lower bound in a universal form.
A universal upper bound for global weak solutions was given in \cite{S3},
and was later recorded and used in \cite{Po}. A local version for smooth 
solutions was obtained by the fourth author, and presented in his course at the
Courant Institute in 1989-90. This result was written up and applied
by Y.-Y. Li \cite{Li}. 

The main subtlety in obtaining an upper bound is the existence of
spherical solutions on $\RR^n$. 
\begin{definition}
For any $\lambda > 0$ and $x_0 \in \RR^n$, the spherical solution
of dilation $\lambda$ and center $x_0$ is given by
\[
u_{\lambda, x_0}(x) =  |x-x_0|^{\q} \big( \cosh ( -\log |x-x_0| +\log \lambda)) \big)^{\q}
=\left( \frac{2\lambda}{1 + \lambda^2 |x-x_0|^2} \right)^{\frac{n-2}{2}}.
\]
\end{definition}
Any one of these functions gives a conformal factor transforming
the flat metric on $\RR^n$ to the pullback of the standard metric
on the sphere by a conformal transformation. The danger is that
a general solution may be well approximated by one of these,
in particular by a strongly dilated one, in some neighbourhood.
This corresponds to the phenomenon of bubbling.

We present here a local version of the upper bound for solutions defined 
in a region which extend to be weak supersolutions across the complement. 
A positive function $u$ in $L^{\frac{n+2}{n-2}}_{\rm{loc}}(B^n)$ is a weak 
supersolution of (\ref{eq:1.1}) if it satisfies 
$$
\Delta u +{n(n-2)\over 4}u^{\frac{n+2}{n-2}}\leq 0
$$
in the distributional sense on $B^n$. We will be especially interested in
supersolutions which are smooth solutions on an open subset $\Omega$ of $B^n$. 
The
complement $\Lambda =B^n\backslash \Omega$ is then relatively closed 
in $B^n$. Note that such a function $u$ is superharmonic, and thus 
can be redefined on a set of measure zero so as to be upper semicontinuous. 
As such, the restriction of $u$ to any compact subset of $B^n$ achieves
its infimum.
%, and in particular, $u$ satisfies
%$$
%\inf_{B(0,3/4)}u=\inf_{\partial B(0,3/4)}u\geq c\int_{B(0,7/8)}u\ dx
%$$
%for some constant $c$ depending only on $n$.

We remark that in the main application, $u$ will be a weak solution. The
reason for considering supersolutions is that it was shown in 
(\cite{SY},Theorem 5.1) that any solution $u$ which defines a complete metric
on a region $\Omega$ in the sphere extends as a weak supersolution to the
full sphere. Moreover, any solution $u$ on $\Omega$ which lies in
$L^{\frac{n+2}{n-2}}_{\rm{loc}}(B^n)$ and tends locally uniformly to infinity
on approach to $\Lambda$ extends as a weak supersolution on $B^n$. This
can be seen by observing that for any large constant $L$, the function
$u_L=\min\{u,L\}$ is a weak solution of the inequality
$$ \Delta u_L \leq -{n(n-2)\over 4}u^{\frac{n+2}{n-2}}\chi_{\{u\leq L\}}
$$
where $\chi_E$ denotes the characteristic function of a set $E$.
A simple application of the dominated convergence theorem then shows that
$u$ defines a weak supersolution.

The proof of the upper bound estimate for a solution $u$ in $B^n$
which is singular only at the origin will require that $u$ extends 
to a weak solution of
(\ref{eq:1.1}) on all of $B^n$. This follows from a more
general result. In fact, any solution of (\ref{eq:1.1}) defined 
outside a  `sufficiently thin' set automatically extends as a weak solution.
To make this precise, we let $\Lambda \subseteq B^n$ be a relatively closed 
set of
Lebesgue measure zero. We will call $\Lambda$ {\bf thin} if there is a 
sequence $\phi_i$ of smooth functions on $B^n$ with values in $[0,1]$ such
that $\phi_i \equiv 0$ in a neighborhood of $\Lambda$, $\lim \phi_i=1$ on
$B^n \setminus \Lambda$, and 
$$ \lim_{i\rightarrow \infty} \int_{B(0,r)}(|\Delta \phi_i|^{(n+2)/4}
+|\nabla \phi_i|^{(n+2)/2})dx=0
$$
for any $r<1$. One sees that a smooth submanifold of dimension less than
$\frac{n-2}{2}$ is thin by choosing $\phi_i$ to be a function of distance
to the submanifold; in particular, a point is thin in any dimension.
\begin{lemma}
Suppose $u$ is a solution of (\ref{eq:1.1}) defined on $B^n\backslash \Lambda$
where $\Lambda$ is a relatively closed thin set. Then $u$ lies locally in 
$L^{(n+2)/(n-2)}(B^n)$, and defines a weak solution of (\ref{eq:1.1}) on $B^n$.
\label{lem:3.1}
\end{lemma}

To prove this we first show that $u$ is locally in $L^{\frac{n+2}{n-2}}(B^n)$.
Let $\zeta$ be a smooth nonnegative function with compact support in 
$B^n\setminus \Lambda$, 
and multiply (\ref{eq:1.1}) by $\zeta^{(n+2)/2}$ and integrate
by parts to obtain
$$ \frac{n(n-2)}{4} \int \zeta^{(n+2)/2}u^{(n+2)/(n-2)}dx=
-\int u\Delta(\zeta^{(n+2)/2}) dx.
$$
This implies the bound 
$$\int (\zeta^{(n-2)/2}u)^{(n+2)/(n-2)}dx\leq 
c\int \zeta^{(n-2)/2}u(\zeta|\Delta \zeta|+|\nabla \zeta|^2)dx.
$$
An application of H\"older's inequality then implies
$$ \int (\zeta^{(n-2)/2}u)^{(n+2)/(n-2)}dx\leq
c\int (\zeta |\Delta \zeta|+|\nabla \zeta|^2)^{(n+2)/4}. 
$$
Now we choose $\zeta=\psi \phi_i$ where $\psi$ is a function which is
equal to one on $B(0,r)$ for some $r<1$, and equal to zero outside $B(0,r_1)$
for some $r_1\in (r,1)$. We then have by an easy estimate 
$$(\zeta |\Delta \zeta|+|\nabla \zeta|^2)^{(n+2)/4}\leq 
c(1+|\Delta \phi_i|^{(n+2)/4}+|\nabla \phi_i|^{(n+2)/2})
$$ 
with a constant $c$ depending on $r,r_1$. Letting
$i$ tend to infinity we then conclude that $u$ is in 
$L^{\frac{n+2}{n-2}}(B(0,r))$ for any $r<1$.

To complete the proof of the lemma, we now let $\zeta$ be any chosen smooth
compactly supported function in $B^n$, and we multiply (\ref{eq:1.1}) by
$\zeta \phi_i$ and integrate by parts to obtain
$$ \int_{B^n}(u\Delta(\zeta \phi_i)+
\frac{n(n-2)}{2}\zeta \phi_iu^{(n+2)/(n-2)})dx=0.
$$
We expand 
$$ \Delta(\zeta \phi_i)=\zeta\Delta(\phi_i)+
2\langle\nabla\zeta,\nabla\phi_i\rangle+\phi_i\Delta\zeta.
$$
By H\"older's inequality 
$$|\int u \zeta\Delta(\phi_i)|\leq \|u\zeta\|_{L^{(n+2)/(n-2)}}
\|\Delta(\phi_i)\|_{L^{(n+2)/4}({\rm spt}\zeta)}, 
$$
and this term tends to zero as $i$ tends to infinity. Similarly,
$$ |\int u \langle\nabla\zeta,\nabla\phi_i\rangle|\leq 
\|u|\nabla\zeta|\|_{L^{(n+2)/(n-2)}}\|\nabla(\phi_i)\|_{L^{(n+2)/4}({\rm spt}
\zeta)}
$$
which also goes to zero. Therefore we may apply the bounded convergence
theorem to let $i$ tend to infinity and conclude
$$ \int_{B^n}(u\Delta(\zeta)+
\frac{n(n-2)}{2}\zeta u^{(n+2)/(n-2)})dx=0.
$$
This completes the proof of Lemma~\ref{lem:3.1}.  

In the following proposition we let $d(x,\Lambda)$ denote the distance from
$x$ to $\Lambda$ for $x\in B^n$.  

\begin{proposition}
Let $u$ be a positive supersolution of (\ref{eq:1.1}) which is a smooth
solution on an
open set $\Omega \subset B^n$. Then there exists a constant $c$ depending 
only on $n$ such that 
$$
u(x) \leq c\,d(x,\Lambda )^{\frac{2-n}{2}}\left(\inf_{\partial 
B(0,3/4)}u\right)^{-1}
$$
for all $x\in \Omega \cap B(0,1/2)$.
\label{pr:3.1}
\end{proposition}

Let $\Lambda_{\epsilon}$ denote the neighborhood of radius $\epsilon$ about
$\Lambda$. For $x\in B(0,5/8)$, set 
$U_{\epsilon}(x)=d_{\epsilon}(x)^{\frac{n-2}{2}}u(x)$
where $d_{\epsilon}(x)=\min\{d(x,\Lambda_{\epsilon}),\frac{5}{8} -|x|\}$. 
In order to prove the proposition it suffices to show that
\begin{equation}
\sup_{B(0,5/8)} U_{\epsilon}\leq c\,\left(\inf_{\partial B(0,3/4)}u
\right)^{-1}
\label{eq:3.1}
\end{equation}
for some constant $c$ not depending on $\epsilon$. 

The first step in proving this is to show that given any constants $R$
and $\delta \in (0,1/2]$, there is a constant $c_1$ depending only on $n$, 
$R$, and $\delta$ such that if 
$M_{\epsilon}\equiv \sup_{B(0,5/8)} U_{\epsilon}\geq c_1$, then 
the solution $u$ ``differs by at most $\delta$'' from a spherical piece 
``of size $R$'' which is 
concentrated near any chosen maximum point $x_0\in B(0,5/8)\backslash 
\Lambda_{\epsilon} $ of $U_{\epsilon}$. Precisely we mean that,
if $x_0$ is a maximum point of $U_\epsilon$, and we set 
$\lambda=u(x_0)^{2/(2-n)}$, then the rescaled function 
$w_\lambda (x) \equiv \lambda^{(n-2)/2}u(x_0 + \lambda x)$ satisfies
\begin{equation}
 \|u_{\mu,y_0}-w_\lambda\|_{C^2(B(0,R))}< \delta
\label{eq:3.2}
\end{equation}
where $u_{\mu,y_0}$ is the standard spherical solution defined above
with $|y_0|\leq c_2$ and $1 /c_2 \leq\mu\leq c_2$ for a constant $c_2$
depending only on $n$. 
To prove this assertion, we use an indirect blow-up 
argument. Assume that there is a sequence $u_i$ with corresponding
values $M_{i,\epsilon}$ tending to infinity. Let $x_i$ be a maximum point of
$U_{i,\epsilon}$, and let $w_i(x)=\lambda_i ^{(n-2)/2}u_i(x_i + \lambda_i x)$ 
where $\lambda_i=u_i(x_i)^{2/(2-n)}$. Since $x_i$ is a maximum point of 
$U_{i,\epsilon}$,
we have $U_{i,\epsilon}(x)\leq U_{i,\epsilon}(x_i)$ for $x\in B(x_i,r_i)$ 
where $r_i=\ha d_{i,\epsilon}(x_i)$. Since 
$d_{i,\epsilon}(x)\geq \ha d_{i,\epsilon}(x_i)$ for $x\in B(x_i,r_i)$ we have
$u_i(x)\leq 2^{(n-2)/2}u_i(x_i)$ for $x\in B(x_i,r_i)$. We then see that $w_i$ 
is a solution of (\ref{eq:1.1}) which is regular and bounded by $2^{(n-2)/2}$
 in $B(0,R_i)$
where $R_i=r_i/\lambda_i =\ha M_i^{2/(n-2)}$. Since $R_i$ tends to infinity,
we see in a standard way that a subsequence of the $w_i$ converges in $C^2$
norm on compact subsets to a global positive solution $w$ of (\ref{eq:1.1}). By
the theorem of \cite{GNN} this must be a standard spherical solution. 

We now argue by contradiction. If the result were not true, there would 
exist a sequence  $u_i$ with corresponding values $M_{i,\epsilon}$ tending 
to infinity such that (\ref{eq:3.2}) does not hold for any $\mu>0$ and any 
$y_0\in {\RR}^n$. This clearly contradicts the previous analysis.

Once (\ref{eq:3.2}) is proved, the bounds on $y_0$ and $\mu$ follow from 
the conditions $w_\lambda \leq 2^{(n-2)/2}$ in $B(0, R')$, where $R'=\ha 
d_\e(x_0)/\lambda \geq  \ha c_1^{2/(n-2)}$ and $w_\lambda (0)=1$. Indeed, 
increasing $c_1$ if necessary, we may always assume that $c_1 \geq 2^{n-2}$ 
so that $R'\geq \ha c_1^{2/(n-2)}\geq 2$.  It follows from  $w_\lambda \leq 
2^{(n-2)/2}$ and $\delta \leq 1/2$ that $u_{\mu,y_0}(y) \leq 2^{(n-2)/2} +
\delta \leq 2^{n/2}$ for all  $y \in B(0, 2) \subset B(0, R')$ and it follows 
from  $w_\lambda (0)=1$ that $u_{\mu,y_0} (0)\geq  1- \delta \geq 1/2$.
 Now assume that $y_0 \in B(0,2)$, we may then take $y=y_0$ in the 
inequality $u_{\mu,y_0}(y) \leq 2^{n/2}$ to obtain $2 \mu \leq 2^{n/(n-2)}$ 
which gives a bound on $\mu$. Finally, if $y_0 \notin B(0,2)$, there exists 
$y_1 \in B(0,2)$ such that $|y_1-y_0|^2 = |y_0|^2 -1$ and it then follows 
from the inequality  $u_{\mu,y_0}(y_1) \leq 2^{n/2}$ that
\[
2^{-2/(n-2)} \mu \leq 1 + \mu^2 |y_0|^2 -\mu^2
\]
and from the inequality $u_{\mu,y_0} (0)\geq 1/2$ that
\[
1 +\mu^2 |y_0|^2 \leq 2^{n/(n-2)} \mu.
\]
These inequalities imply first that  $ 2^{-n/(n-2)} \leq \mu \leq 2^{n/( n-2)}
$ and then $|y_0|\leq 2^{n/(n-2)}$.

Notice that the conclusion (\ref{eq:3.2}) can be replaced by the simpler
inequality
\begin{equation}
 \|u_{1,0}-w_{\lambda} \|_{C^2(B(0,R))}< \delta
\label{eq:3.3}
\end{equation}
if we allow ourselves to shift the center point $x_0$ to a new point $x_1$
which is within distance $2c_2\lambda$ where $c_2=c_2(n)$ is the constant 
above, and which is a local maximum point of $u$. Thus we now are taking 
$w_{\lambda} (x) \equiv \lambda^{(n-2)/2}u(x_1 +
 \lambda x)$ with $\lambda =u(x_1)^{2/(2-n)}$. To see this 
we observe from the blow-up argument above that $w_i$ will have a 
absolute maximum point near $y_0$ for $i$
sufficiently large, since $|y_0|$ is bounded and since $\mu$ is bounded from below and from above. Shifting center to this point, and rescaling so that
the maximum value is one then produces $u_{1,0}$ as the limit of the sequence.
Note that the distance $2c_2\lambda_i$ is an arbitrarily small constant
times $r_i$ since $\lambda_i=2M_i^{2/(2-n)}r_i$, and thus $x_1$ can be used
in place of $x_0$ above and we can assume, increasing $c_1$ if  necessary, that  $2c_2\lambda_i< 1/16$.

To complete the proof of (\ref{eq:3.1}), we first note that 
$\inf_{\partial B(0,3/4)}u$ is bounded above since one can see easily from
the fact that $u$ is a supersolution (as in the proof of Lemma~\ref{lem:3.1}) 
that $u$ is locally bounded in $L^{(n+2)/(n-2)}(B^n)$. Thus 
if the left hand side is bounded above by the constant
$c_1$ (for a fixed chosen $R$ and $\delta$) of the previous paragraph, then 
we are done. Therefore we may assume
that this is not the case, and that $u$ satisfies the above bounds with
$\delta$ chosen as small as desired, and $R$ chosen as large as we like.
In this case, we are going to prove that (\ref{eq:3.1}) still holds.
We now introduce the function $v(t,\theta)$ on the cylinder 
$\RR\times S^{n-1}$ which corresponds to the function $w=w_\lambda$ 
with $\lambda=u(x_1)^{2/(2-n)}$. Thus we define
$v(t,\theta)=|x|^{(n-2)/2}w(x)$, where $t= - \log |x| $ and $\theta =x/|x|$.
It is important to note in the following argument that because the point
$x_1$ lies in the ball of radius $5/8 + 1/16 = 11/16$ about the origin, its distance to
$\partial B(0,3/4)$ is at least $1/16$; upon rescaling, we see that the
ball $B(0,\frac{1}{16}
\lambda^{-1})$ lies in the domain of $w$, and corresponds to
a portion of $B(0,3/4)$ in the domain of $u$. 
From (\ref{eq:3.3}), we see that $v$ is close in $C^2$ norm on 
$[-\log R,\infty)\times S^{n-1}$ to the function 
$v_0(t,\theta)=(\cosh t)^{(2-n)/2}$. We now {\it fix} $R=e^2$ so that
$-\log R=-2$.
It follows that $\del_t v ( - 1,\theta ) > 0$ for all $\theta \in S^{n-1}$. 
Then apply the Alexandrov technique to $v$ on the region
$[-\log(\frac{1}{16}\lambda^{-1}),\infty)\times S^{n-1}$, reflecting across
the spheres $\{t_1\} \times S^{n-1}$, starting with $t_1$ very positive, 
and continuing as far as possible. Because $\del_t v (-1,\theta )>0$, 
this procedure must end before $t_1$ reaches $-1$.
We will show below that because $v$ is a supersolution, and the reflected 
function $v^*$ is a solution, there can be no interior contact point with
$v^*\leq v$. Since $v$ is a regular solution near $t=t_1$ for $t_1\geq -1$, 
the Hopf boundary point lemma implies that $\del_t v(t_1,\theta ) < 0$. 
Thus the reflection argument must end because a contact between 
$v^*$ and $v$ occurs on the boundary. In order for this to happen we must
have 
$$ \inf \{v(-\log(\frac{1}{16}\lambda^{-1}),\theta):\ \theta\in S^{n-1}\}<
\sup \{v(t,\theta):\ t>\log(\frac{1}{16}\lambda^{-1})-2,\ \theta\in S^{n-1}\}.
$$
Recalling the definitions of $v$ and $\lambda$, we see that
$$ \inf \{v(-\log(\frac{1}{16}\lambda^{-1}),\theta):\ \theta\in S^{n-1}\}=
16^{(2-n)/2}\inf_{\partial B(x_1,1/16)}u\geq 16^{(2-n)/2}\inf_{\partial B(0,3/4)}u 
$$
where the inequality holds because $u$ is superharmonic. Finally, 
for $\delta$ small and {\it fixed}, we have 
$v(t,\theta)\leq 2v_0(t)$ for $t\geq 0$, and therefore we have 
$$ \inf_{\partial B(0,3/4)}u\leq c\sup_{[\log(\frac{1}{16}\lambda^{-1})-2,
\infty)}v_0\leq c\lambda^{(n-2)/2}=c\,u(x_1)^{-1} 
$$
where we have used the definition of $v_0$ and $\lambda$.
This implies (\ref{eq:3.1}) since
$U_{\epsilon}$ is bounded above by a constant times $u$.

In this argument we needed to know that if $v^*\leq v$ in a connected 
open set, then either $v^*<v$ or $v^*\equiv v$. To see this, observe that 
$L(v-v^*)\leq (v^*)^{(n+2)/(n-2)}-v^{(n+2)/(n-2)}\leq 0$ weakly. 
The standard mean value inequality then gives us the desired conclusion. 
This completes the proof of the proposition.  

For the proof of the lower bound, we require the following result.
\begin{lemma}
If $u$ is a positive solution of (\ref{eq:1.1}) defined in $B^n\backslash 
\{0\}$ with $\lim_{x\rightarrow 0}|x|^{(n-2)/2}u(x)=0$, then 
$u$ extends as a smooth solution to all of $B^n$.
\label{lem:3.2}
\end{lemma}

To prove this, first note that by Lemma~\ref{lem:3.1} $u$ extends 
as a weak solution to $B^n$. If we could show that $u\in L^p$ for 
some $p>2n/(n-2)$, then standard results would
imply first that $u$ is bounded near the origin, and then by 
linear elliptic theory that it extends smoothly across $0$. 
Let $v(t,\theta )$ be the corresponding cylindrical solution, and observe that
the hypothesis implies that $v(t,\theta )$ tends to $0$ uniformly 
as $t$ tends to $\infty$. Therefore (\ref{eq:1.1}) implies that
$\Delta v\geq \beta v$ for some $\beta>0$ and $t\geq t_0$ for 
sufficiently large $t_0$. Now consider the function 
$w=ce^{-\sqrt{\beta}t}+\epsilon e^{\sqrt{\beta}t}$ where $c$ is
chosen large enough that $ce^{-\sqrt{\beta}t_0}>v(t_0,\theta)$ 
for all $\theta \in S^{n-1}$. Note that $w$ is a solution of 
$\Delta w=\beta w$. We will use the maximum principle on 
$[t_0, T]\times S^{n-1}$ for any sufficiently large $T$ to obtain the
necessary decay estimate on $v$. If $T$ is chosen sufficiently positive, then 
we have $v(T,\theta)<w(T,\theta)$ since $v$ tends to zero as $t$ tends to 
$\infty$ while $w$ is growing at an exponential rate. 
By the maximum principle it follows that for all $t_0\leq t\leq T$, all 
$\theta$, and any chosen $\epsilon >0$ we have
$$ 
v(t,\theta) \leq ce^{-\sqrt{\beta}t}+\epsilon e^{\sqrt{\beta}t}.
$$
Since $T$ is arbitarily large, this inequality holds for all $t\geq t_0$, and
since $c$ is independent of $\e$, we may let $\e$ go to zero.
We have shown that $e^{\sqrt{\beta}t}v(t,\theta) \le c$ for
$t\geq t_0$. Writing this in terms of $u$ we have shown that
$u(x)\leq c|x|^q$ for $q=(2-n)/2+\sqrt{\beta}$. 
This implies that $u\in L^p$ for some $p>2n/(n-2)$ as required. 
This completes the proof of Lemma~\ref{lem:3.2}.

We now are in a position to prove Theorem~\ref{th:3.1}. 
The upper bound follows from Lemma~\ref{lem:3.1} and 
Proposition~\ref{pr:3.1}, as already noted. To prove the 
lower bound, it is most convenient to use the cylindrical
background, and so consider the function $v(t,\theta)$.
We first observe that the upper bound $v \le C_2$ implies 
a Harnack inequality of the form 
\begin{equation}
\sup\{v(t,\theta):\ T \le t \le T + T_0\}\leq 
c\,\inf \{v(t,\theta):\ T \le t \le T+T_0\}
\label{eq:ll}
\end{equation}
for all $T > 2$, and where the constant $c$ depends
on $T_0$ and $C_2$, but not on $T$. 
Define $\bar{v}(t)=\int_{S^{n-1}}v(t,\theta )\ d\theta$. 
From the Harnack inequality we see that $v$ and $\bar{v}$ 
are bounded in ratio. Suppose then that
the singularity at $0$ is not removable. Applying Lemma~\ref{lem:3.2} we see 
that either the lower bound holds or else there is a sequence
$t_i$ tending to $\infty$ such that $\bar{v}(t_i)$ tends to $0$. 

Using (\ref{eq:ll}), the Harnack inequality, standard elliptic estimates as well as the fact that $v$ is a solution of (\ref{eq:1.1}), we can estimate the radial Poho\^zaev invariant ${\cal P}_{\rm ra}(v)$,
as defined in \S 2 (cf. also \S 6), at any one of the $t_i$.
Using the fact that it is independent of $i$, we see
that the invariant must vanish. 
Now set $w_i(t,\theta)=v(t+t_i,\theta)/\bar{v}(t_i)$.
Some subsequence of the $w_i$ converges in $C^2$ on compact subsets of
$\RR \times S^{n-1}$ to a positive solution of $Lw=0$. The associated 
function $h(x)$ defined on $\RR ^n\setminus \{0\}$ is then a positive 
harmonic function, and can therefore be written $a|x|^{n-2}+b$ for some 
$a,b\geq 0$. This is the same as the condition that $w(t,\theta)
=ae^{(n-2)/2t}+be^{-(n-2)/2t}$.  Since $\bar{w}$ has a critical 
point at $t=0$, it follows that $a = b$, and both are positive. 
Then, a direct calculation shows that 
\[
\lim_{i\rightarrow +\infty}
\int_{S^{n-1}} (\frac{1}{2} (\del_t w_i)^2 -\frac{1}{2}  
|\nabla_\theta w_i |^2 
- \frac{(n-2)^2}{8} w_i^2 + \bar{v}(t_i)^{\frac{4}{n-2}}
\frac{(n-2)^2}{8} w_i^\pp )d\theta
\]
\[
= \int_{S^{n-1}} \frac{1}{2} (\del_t w)^2 
- \frac{(n-2)^2}{8} w^2  d\theta =- \omega_n  \frac{(n-2)^2}{2} ab \neq 0
\]
But this is a contradiction, for this is also a limit of rescalings of 
the Poho\^zaev invariant for $v$, and hence must be zero.
Indeed
\[
0=\bar{v}(t_i)^{-2}  {\cal P}_{\rm ra}(v)= 
\int_{S^{n-1}} (\frac{1}{2} (\del_t w_i)^2 -\frac{1}{2}  
|\nabla_\theta w_i |^2 - \frac{(n-2)^2}{8} w_i^2 + 
\bar{v}(t_i)^{\frac{4}{n-2}}
\frac{(n-2)^2}{8} w_i^\pp) d\theta \neq 0
\]
for $i$ sufficiently large.  This establishes the lower bound. 

To show that ${\cal P}_{\rm ra}(v) \neq 0$, we observe that by the 
upper and lower bounds, we can choose a sequence $t'_i$ tending to $\infty$ 
so that the corresponding translated solutions $t\rightarrow v(t+t'_i, \theta)$
converge in $C^2$ norm on compact subsets of $\RR \times S^{n-1}$ to a solution
satisfying the same bounds and defined on all of $\RR \times S^{n-1}$. By 
Proposition~\ref{pr:2.1} such a solution must be independent of $\theta$,
and hence must have nonzero radial Poho\^zaev invariant. Since this invariant
must be the same as ${\cal P}_{\rm ra}(v)$ because of the $C^2$ convergence,
we have shown ${\cal P}_{\rm ra}(v) \neq 0$
as desired. The proof of Theorem~\ref{th:3.1} is complete.

\section{Simple convergence to a radial solution}

Our aim in this section is to prove 
\begin{proposition} 
Assume that $u$ is a nonnegative smooth solution of (\ref{eq:1.1}) 
defined in the punctured unit ball $B^n\backslash \{0\}$ with a nonremovable
singularity at the origin. For $t>0$ and $\theta \in S^{n-1}$, let us 
define $v(t, \theta)$ so that $u(x)=|x|^{\frac{2-n}{2}} v(-\log |x|, x/|x|)$.
 Then, there exists 
a Fowler parameter $\e \in (0,\e_0]$, constants $T \in [0,T_\e)$ and $C > 0$,
and some exponent $\alpha > 0$ such that
\[|v(t,\theta) - v_\e(t+T)| \le C e^{-\alpha t}, \quad \mbox{ for }\ t \ge 0.\] 
\label{pr:simasym}
\end{proposition}

We first show that any sequence of translates of $v$ has a subsequence 
converging to one of the $v_\e$, then that any angular derivative of 
$v$ converges to zero. After this, a somewhat delicate rescaling argument
due originally to Leon Simon, in a different context, 
gives the final convergence. 

By the results of the last section we know that $0 < C_1 \le v(t,\theta) 
\le C_2$ for all $t \ge 0$. By standard elliptic estimates, we also
get the uniform boundedness of any derivative $|\del_t^j 
\del_\theta^\gamma v(t,\theta)| \le C_{j,\gamma}$ for all $t \ge 0$. 

Let $\{\tau_j\}$ be any sequence
of numbers converging to $\infty$, and define $v_j(t,\theta) =
v(t + \tau_j,\theta)$. Then $v_j$ is defined on $[-\tau_j,\infty)
\times S^{n-1}$ and satisfies (\ref{eq:2.1}) there. Using
the uniform bounds on any derivative of $v_j$, we may choose a 
subsequence of the $v_j$ converging in the ${\cal C}^\infty$ topology on any
compact subset of $\RR \times S^{n-1}$. The limit function, $v_{\infty}$, 
still satisfies (\ref{eq:2.1}), is nonvanishing because of the lower
bound for $v$ and is defined on the whole cylinder. By 
Proposition~\ref{pr:2.1}, the only functions with these properties are
the translated Fowler solutions; hence we deduce that 
for some $\e$ and $T$, $v_{\infty}(t,\theta) = v_\e(t + T)$. 

The radial Poho\^{z}aev invariant of $v$ equals that of any one of 
the $v_j$, hence also equals that of $v_\e(t+T)$. Hence the limiting 
necksize $\e$ is independent of the original sequence of numbers 
$\tau_j$ and of the subsequence. 

The fact that any sequence of translates of $v$ has a subsequence
converging to a $\theta$-independent solution is very strong. 
It implies immediately that any angular derivative
$\del_{\theta} v$ tends to zero uniformly. For if this were false,
then there would be a sequence of points $(\tau_j,\theta_j)$ with
$\tau_j \rightarrow \infty$ for which $|\del_{\th} v(\tau_j,\theta_j)| 
\ge C > 0$. But then, translating back by $\tau_j$ and rotating
$\theta_j$ to some fixed point $\theta_0 \in S^{n-1}$ to get
a new sequence of solutions $v_j$, we can
again extract some subsequence converging to a radial function.
But this contradicts the positive lower bound on 
$|\del_{\theta} v_j(0,\theta_0)|$.

Next, let $X$ be any infinitesimal rotation on $S^{n-1}$, which we
denote for simplicity by $\del_\theta$; 
applying it to (\ref{eq:2.1}) and using that 
it commutes with $\Delta_\th$, we see that $\del_\theta v$ is a Jacobi field 
for the Jacobi operator at $v$, i.e. the linearization of (\ref{eq:2.1}) 
at $v$. By the discussion above, we know that $\del_\theta v = o(1)$, but 
since we know little about this Jacobi operator, we cannot deduce better
decay directly. However, for some $\tau_j \rightarrow \infty$ 
consider the corresponding sequence of translates $\del_\th v_j$ and 
let $A_j = \sup_{t \ge 0} |\del_\th v_j|$. If this supremum
is attained at some point $(s_j,\th_j)$ then $s_j$ must stay
bounded. Otherwise we could translate back further by $s_j$ to
obtain, after passing to a subsequence, a Jacobi field $\phi$ 
for the Jacobi operator at $v_\e(t+T)$, for some $T$, defined
on all of $\RR$. By construction, $\phi$ is bounded and not
identically zero, but on the other hand, since it clearly has no zero 
eigencomponent relative to $\Delta_\th$, the results above show
that it must increase exponentially either as $t \rightarrow +\infty$
or as $t \rightarrow -\infty$. This contradiction shows that
$s_j$ stays bounded. Now fix a positive integer $N$ such that
$NT_\e > \sup s_j$ and define $I_N \equiv \{0 \le t \le NT_\e\}$; here,
of course, $T_\e$ is the period of $v_\e$. Then $A_j = \sup_{I_N} |\del_\th
v_j|$, and we obtain again, after passing to a subsequence, a Jacobi
field $\phi$ for the Jacobi operator at some $v_\e(t+T)$ which is
bounded for $t \ge 0$ and attains its supremum of $1$ in $I_N$. 

By the results on this Jacobi operator stated earlier,
$\phi$ must decay at least like $e^{-t}$. Actually, something
slightly stronger is true. Although the limit $\phi$ is not necessarily
unique, there exists a constant $c$, independent of all other choices, 
such that $|\phi(t)| \leq c e^{-t}$ for $t \geq 0$.
To see this, write $\phi = \sum_{j=1}^n c_j \psi_{\e,j}^+ + \phi'' = \phi' + 
\phi''$, where $\phi'' \in \mbox{Span} \{ h_j(t) \chi_j(\theta): j 
\geq n+1 \} \equiv E''$. Since $|\phi| \le 1$ for $t \ge 0$, the $c_j$, 
$j=1, \ldots , n$, are absolutely bounded because 
they are given by the (normalized) 
inner products of $\phi$ with $\chi_j(\theta)$.  
Hence $|\sum_{j=1}^n c_j \psi_{\e,j}^+| \le c e^{-t}$ with $c = c(n)$. 
Next, we claim that $|\phi''| \leq ce^{-t}$ for some $c$ independent
of all other choices, and that this is a consequence of the fact
that $\phi''$ is absolutely bounded for $t \ge 0$ and 
$\phi'' \in E''$. We argue by contradiction and assume that 
there exists a sequence of Jacobi fields $\phi_i'' \in E''$ such that 
$A_i'' = \sup_{t\geq 0} e^{t}| \phi_i''(t,\theta)|$ tends to infinity. 
If this supremum is achieved at a point $(t_i,\theta_i)$, then define
$\tilde{\phi}_i(t,\theta) = e^{t_i}A_i^{-1}\phi_i''(t + t_i,\theta)$; these
are Jacobi fields on $[-t_i,\infty)$ which are bounded by $e^{-t}$, 
and which take the value $1$ at $(0,\theta_i)$. 
The $t_i$ clearly cannot be bounded, because the Jacobi equation
is uniformly elliptic on any finite piece of the cylinder. 
Passing to a subsequence, if necessary, we obtain a Jacobi field 
$\tilde{\phi}$ defined and satisfying 
$|\tilde{\phi}| \le e^{-t}$ on the whole cylinder. Since we
also know that $\tilde{\phi} \in E''$, it must grow faster
than $e^{-t}$ either at $+\infty$ or $-\infty$, which is
a contradiction. This proves the claim. 

At this stage, it is now possible to show that $\del_\th v$ decays 
exponentially.
Although we do not, strictly speaking, need this result, we sketch
the proof anyway as a warm-up to the arguments later. Let 
$J_N \equiv \{NT_\e \le t \le 2NT_\e\}$, and denote by $v_\tau(t,\th)$
the translate $v(t+\tau,\th)$. We claim that 
if $N$ is chosen large enough then
\[ 
\forall \tau \ge 0, \qquad \sup_{J_N} |\del_\th v_\tau| \le \ha
\sup_{I_N} |\del_\th v_\tau|.
\]
As usual, this is proved by contradiction. If it were to fail, there
would exist some sequence $\tau_j \rightarrow \infty$ such that
$B_j > \ha A_j$, where
\[
A_j \equiv \sup_{I_N} |\del_\th v_j|, \qquad B_j \equiv
\sup_{J_N} |\del_\th v_j|.
\]
However, since $A_j^{-1}\del_\th v_j$ converges to the Jacobi
field $\phi$ uniformly on compact sets, we see that $\del_\th v_j = 
A_j \phi + o(A_j)$ on $I_N \cup J_N$. Since $|\phi| \le Ce^{-t}$, 
we then get that $B_j \le C A_j e^{-NT_\e} + \frac{1}{4}A_j$ for
$j$ sufficiently large. However, as we have seen,  the constant
$C$ is bounded independently of other choices since $\phi$ is normalized
to be less than $1$ on $I_N$, so if $N$ is sufficiently large, 
$Ce^{-NT_\e} \le \frac{1}{4}$, and we obtain a contradiction.
This proves the assertion.  It is now trivial to deduce 
that $\del_\th v$ decays at some exponential rate. 

We may use the sequence of translates $v_j$ to obtain a Jacobi field 
in another way. Suppose that we have chosen a subsequence converging to some 
$v_\e(t+T)$. Define
\[
w_j(t,\theta) = v_j(t,\theta) - v_\e(t+T), \qquad \alpha_j = 
\sup_{I_N} |w_j|, \qquad \mbox{\rm and } \quad 
\phi_j = \alpha_j^{-1} w_j.
\]
Then it is easy to check that $\phi_j$ converges to a solution 
$\phi$ of the Jacobi operator at $v_\e(t+T)$. 
This Jacobi field is bounded on $t \ge 0$. To see this we
examine its eigencomponents with respect to $\Lap_\theta$.
First consider the sum of eigencomponents over all nonzero
eigenvalues, which we call $\tilde{\phi}$. To show that
$\tilde{\phi}$ is bounded, and hence exponentially 
decaying, it is sufficient to show that $\del_\theta \phi$
is also bounded for $t \ge 0$. We may assume that $\tilde{\phi}$
is not identically zero, for otherwise this case is trivial. 
This function arises as the limit of a subsequence of the sequence 
$\alpha_j^{-1} \del_\theta v_j$. If we knew that $\alpha_j$
were commensurate with $A_j$, the supremum of $\del_\th v_j$
on $I_N$, i.e. if $C_1 A_j \le \alpha_j \le C_2 A_j$ for some constants
$C_1, C_2 > 0$, then we could appeal to our earlier argument.  However,
this must be the case, for otherwise $\alpha_j^{-1} \del_\th v_j$
would either blow up or tend to zero on $I_N$, whereas we know
that it converges to the nontrivial function $\del_\th \phi$. As for
the remaining eigencomponent, with eigenvalue zero, we know that
\[
\phi(t,\theta) = a^+ \psi_{\e,0}^+ + a^- \psi_{\e,0}^- + \tilde{\phi},
\]
where $\tilde{\phi}$ decreases exponentially. We claim that $a^- = 0$. 
Intuitively this is true because this coefficient corresponds to an 
infinitesimal change of Fowler parameter. Rigorously, simply write 
$v_j = v_\e + \alpha_j \phi + o(\alpha_j)$ on $I_N$, say; if $a^-$ were not to 
vanish, then the radial Poho\^{z}aev invariants of the two sides would 
differ, which is impossible. Thus we have shown that $\phi$ must
be bounded for $t \ge 0$. 

We are finally in a position to show that the difference between
$v(t,\th)$ and some $v_\e(t+T)$ converges to zero uniformly.
The subtlety here arises because the displacement $T$ cannot be detected
by any of the Poho\v{z}aev integrals. We use again the same sort of 
argument as above, which is due to L. Simon. 
Define $v_\tau(t,\theta) = v(t+\tau,\theta)$ and set
$w_\tau(t,\theta) = v_\tau(t,\theta) - v_\e(t)$. (Since we do not
know the correct translation parameter $T$ beforehand, we choose
one arbitrarily, say $T=0$.)  The idea is to establish
an improvement of approximation in the following sense. 
Use the interval $I_N$ as before and fix a constant $B > 0$, and let 
$\eta(\tau) = \sup_{I_N} |w_\tau|$. Then we claim~: 

If $\tau$ is sufficiently large and $\eta(\tau)$
sufficiently small, then there exists an $s$ with $|s| \le B 
\eta(\tau)$ such that $\eta(\tau + NT_\e + s) \le \frac{1}{2}\eta(\tau)$.

To prove this, suppose that, for fixed values of $N$ and $B$, it fails.
Then there exists some sequence $\tau_j \rightarrow \infty$ such
that $\eta_j \equiv \eta(\tau_j)\rightarrow 0$, but such that for any
$s$ with $|s| \le B\eta_j$ we have $\eta(\tau_j + NT_\e +s) > \frac{1}{2}
\eta_j$. Let $\phi_j = \eta_j^{-1}w_{\tau_j}$. We have shown already
that $\phi_j$ converges in ${\cal C}^\infty$ on compact sets
to a Jacobi field $\phi$ which is bounded for $t \ge 0$. By construction,
$|\phi| \ge 1/2$ on $NT_\e \le t \le 2NT_\e$, so $\phi \not\equiv 0$. 
Expanding
$\phi$ into eigenfunctions again, we get $\phi = a^+\psi_{\e,0}^+ + \tilde{\phi}$,
where $\tilde{\phi}$ is exponentially decreasing. Note that $a^+$, which
we relabel simply $a$, is uniformly bounded, independently of the 
sequence, because $|\phi| \le 1$ on $0 \le t \le NT_\e$. We assume
that the constant $B$ has been chosen larger than this upper bound. 
Because $\psi_{\e,0}^+$  
corresponds to infinitesimal translation, this argument indicates
that we should adjust $v_\e$ by some translation, and to first order,
this adjustment should just be $s_j = -\eta_j a$; by our choice of
$B$, this is less than $B\eta_j$. Thus we estimate
\[
w_{\tau_j + s_j}(t,\theta) = v(t+\tau_j - \eta_j a,\theta) - v_\e(t)
= w_{\tau_j}(t,\theta) - a \eta_j \psi_{\e,0}^+(t) + o(\eta_j).
\]
In particular,
\[
w_{\tau_j + s_j} = \eta_j \tilde{\phi} + o(\eta_j).
\]
But the function $\tilde{\phi}$ decays at a fixed exponential
rate and its supremum on $0 \le t \le NT_\e$ is bounded by one, so that
if $N$ is chosen large enough (but again, independent of our original
sequence $\tau_j$, etc.) then
\[
\sup_{0 \le t \le NT_\e} |w_{\tau_j + s_j + NT_\e}|
= \sup_{NT_\e \le t \le 2NT_\e} |w_{\tau_j + s_j}| 
= \eta_j \sup_{NT_\e \le t \le 2NT_\e} |\tilde{\phi}| + o(\eta_j)
\le \frac{\eta_j}{4}.
\]
But this contradicts the assumption that 
$\eta(\tau_j + NT_\e + s) > \frac{1}{2} \eta(\tau_j)$, and
thus we have proved our claim.

It remains only to show how this approximation improvement property
yields the correct simple asymptotics.
In fact, we must prove that $w_\sigma \rightarrow 0$ for some fixed
translation $\sigma$. We may assume, by starting at some sufficiently
large $t$, that $\eta(0)$ is sufficiently small and that
$B\eta(0) \le \frac{1}{2}NT_\e$. Let $\tau_0 = 0$ and $s_0$ be
determined by the first step of this argument. Define, for any
integer $j > 0$,
\[
\sigma_j = \sum_{i=0}^{j-1} s_i, \qquad \tau_j = \tau_{j-1} + s_{j-1} +
NT_\e.
\]
Then by iteration, $\eta(\tau_j) \le 2^{-j} \eta(0)$, $|s_j| \le 2^{-j-1}
NT_\e$, and hence in particular $\sigma = \lim \sigma_j$ is well defined 
and less than $NT_\e$. To show that $\sigma$ is the correct
translation parameter, for any $t > 0$, let $[t]$ denote the
reduction mod $NT_\e$, so that $t = jNT_\e + [t]$ for some $j \ge 0$. 
Then
\begin{eqnarray*}
w_{\sigma}(t,\theta) = v(t+\sigma,\theta) - v_\e(t)
= \big(v(t+\sigma_j,\theta) - v_\e(t)\big) + \big(v(t+\sigma,\theta)
- v(t+\sigma_j,\theta) \big) \\
= w_{\tau_j}([t],\theta) + O(2^{-j}),
\end{eqnarray*}
where we have used Taylor's theorem and the uniform boundedness
of $\del_t v$ and $\del_t v_\e$. But using the bound on $\eta(\tau_j)$,
we finally conclude that
\[
|w_\sigma(t,\theta)| \le C2^{-j}, \qquad \mbox{\rm or\  equivalently} 
\quad |w_{\sigma}(t,\theta)| \le C' e^{-\frac{ \log 2}{NT_\e} t}.
\]
This completes our proof of Proposition~\ref{pr:simasym}.

\section{Refined asymptotics} 
In this last step of the proof of Theorem~\ref{th:refas}, we improve the 
asymptotics one step further.  We show that for some $a \in \RR^n$, 
$u$ has an expansion of the form 
\begin{equation}
u(x) = |x|^\q\big( v_\e(-\log |x| + T) + 
(a \cdot x)\psi_{\e,1}^+(-\log |x| + T) + O(|x|^\beta) \big), \label{eq:5.1} 
\end{equation}
\[
\mbox{ where } \quad \beta = \min\{2, \gamma_2^+\}. 
\]
On the other hand, we also know that the deformed Fowler solution 
$u_{\e,a}$ has an expansion of this exact form, (\ref{eq:2.25}); 
comparing these two expansions completes the proof of the main theorem.

Using the simple asymptotics already established, write 
\[
u(x) = |x|^\q v(-\log |x|) = |x|^\q\big(v_\e(-\log |x| + T) + 
w(-\log |x|)\big)
\]
where $w(t) \in {\cal C}^{2,\al}_{-\gamma}(\RR^+ \times S^{n-1})$
for some $\gamma >0$. The function $v(t)$ satisfies (\ref{eq:2.1}), 
which we write as $N(v) = 0$. Expanding this in a Taylor series about 
$v_\e(t+T)$ gives
\begin{equation}
\Li_\e w = \frac{n(n-2)}{4}\bigg( (v_\e + w)^\p -
v_\e^\p - \frac{n+2}{n-2}v_\e^{\frac{4}{n-2}} w \bigg) \equiv Q(w).
\label{eq:5.55}
\end{equation}
It is straightforward to check that if $w \in {\cal C}^{0,\al}_{-\gamma}
(\RR^+ \times S^{n-1})$
and $v_\e + w > 0$, then $Q(w) \in {\cal C}^{0,\al}_{-2\gamma}(\RR^+ 
\times S^{n-1})$. 

First assume that $0 < \gamma < \ha$. Then Corollary~\ref{cor:2.rinv} 
gives that $w \in {\cal C}^{2,\al}_{-2\gamma}(\RR^+ \times S^{n-1})$, and 
so $Q(w) \in
{\cal C}^{2,\al}_{-4\gamma}(\RR^+ \times S^{n-1})$. 
Continuing on, we deduce after finitely 
many steps that $w \in {\cal C}^{2,\al}_{-\gamma'}(\RR^+ \times S^{n-1})$
for some $\gamma' \in (1/2, 1)$. Applying Corollary~\ref{cor:2.rinv}
again gives $w \in {\cal C}^{2,\al}_{-\beta}(\RR^+ \times S^{n-1}) 
\oplus E_{\e,1}$, where $\beta = \min\{2\gamma',\gamma_2\})$.
 
The optimal $\beta$ we could expect is $\min \{2,\gamma_2\}$. 
We have proved the expansion (\ref{eq:5.1}). 

\section{The global balancing formula} 

In this section we shall present an application of the refined 
asymptotics theorem. We consider solutions of the
singular Yamabe problem with discrete singular set $\Lambda$ 
in the standard conformal class $[g_0]$ on $S^n$. More specifically, let
$\Lambda = \{p_1, \dots, p_k\} \subset S^n$ be arbitrary, $k \ge 2$.
The solutions we consider are functions $u > 0$ on $\SL$ such that 
$u^{\ppp}g_0$ is a complete metric of constant positive scalar
curvature $n(n-1)$ on $\SL$. More analytically, these functions
satisfy the special case of (\ref{eq:1.3}), 
\begin{equation}
\Lap_{S^n}u - \frac{(n-2)^2}{4}u + \frac{(n-2)^2}{4}u^\p = 0,
\label{eq:10.2}
\end{equation}
with $u > 0$ on $\SL$ and $u$ singular at the points of $\La$. The 
`unmarked moduli space' $\Mk$ consists of all such solutions with $k$ 
singular points, but with $\La$ allowed to vary. (One may also
consider the `marked moduli space' $\ML$, consisting of all such 
solutions with fixed singular set $\La$.)
The most basic question about these moduli spaces, whether they
are nonempty, was answered originally by the fourth author \cite{S1} some
time ago; an alternate construction was recently given in \cite{MP2}. 

The refined asymptotics theorem implies that to any $u \in \Mk$ one can 
associate a set of parameters: the singular points $p_j$, and then
at each $p_j$ the Fowler and translation parameters $\e_j \in (0,\e_0]$ 
and $a_j \in \RR^n$. This yields a map
\begin{equation}
\Phi: \Mk \longrightarrow {\cal X}_k \equiv 
\left( \big(\prod_{j=1}^k S^n\big) \setminus \cup_{i,j}\Delta_{ij}
\times (0,\e_0]^k \times (\RR^n)^k \right) / \Sigma_k,
\label{eq:10.1}
\end{equation}
where $\Delta_{ij}$ is the diagonal of the product of the $i^{\rm th}$
factor and the $j^{\rm th}$ factor and $\Sigma_k$ is the permutation
group on $k$ letters, which acts in the obvious way. It is shown
in \cite{MPo} (following the proof of the analogous fact for $\ML$
in \cite{MPU1}) that $\Mk$ is a real analytic set of dimension $k(n+1)$; 
on the other
hand, $\dim {\cal X}_k = k(2n+1)$. There are actually another $k$
parameters, arising from the translations (what we called $T$ earlier)
at each $p_j$. Thus the correct parameter space is $2k(n+1)$-dimensional.
In current work, D. Pollack shows that this extended
version of $\Phi$ is a Lagrangian immersion which is equivariant under
the obvious action of the conformal group on either side. 

Our purpose here is to use the refined asymptotics theorem
to show that the image of $\Mk$ in ${\cal X}_k$ lies in the
zero set of a collection of real analytic equations obtained by computing
the global Poho\^{z}aev integral for each solution $u$ with
respect to any conformal Killing vector field $X$ on $S^n$. 
Although there are not enough equations here in general to
fully determine the moduli space $\Mk$, the relations
we can obtain here are still interesting in their own right.
Before embarking on this computation, though, we first
define these integral invariants.

\subsection{Poho\^{z}aev invariants} 
We now turn to a discussion of the existence and specific 
form of a family of homological integral invariants of solutions of 
equation (\ref{eq:1.1}). These `Poho\^zaev invariants' were
discovered in their simplest form by Poho\^{z}aev, and applied by him 
to prove the nonexistence of nonnegative smooth solutions of 
(\ref{eq:1.1}) satisfying Dirichlet conditions on smooth,
star-shaped domains. The identity he discovered was later 
put into a natural and general Riemannian setting in \cite{S1}. 
These extended invariants associate to any metric $g$ of constant scalar 
curvature on a manifold $\Omega$, any closed hypersurface $\Sigma 
\subset \Omega$ and any conformal Killing field $X$, a 
real number ${\cal P}(X,\Sigma)$ which depends only on
the the homology class of $\Sigma$ in $H^{n-1}(\Omega)$. 
The invariant associated to $X$, $g$ and $\Sigma$ is
\begin{equation}
{\cal P}(X,g,\Sigma) \equiv \frac{n-2}{2(n-1)}
\int_{\Sigma} T(X,\nu)\,d\sigma,
\label{eq:genpoh}
\end{equation}
where $T(\cdot, \cdot)$ is the trace-free Ricci tensor for the
metric $g$, $\nu$ is the unit normal to $\Sigma$ and $d\sigma$
is the volume form induced by $g$ on $\Sigma$. 

In this subsection we shall compute these invariants 
in our particular case of interest, where $g = u^{\frac{4}{n-2}}\delta$,
$u$ is a solution of (\ref{eq:1.1}) in the ball $B(0,1)$ with an 
isolated singularity at $0$ and $\Sigma = \Sigma_{|\eta|}$ 
is the sphere $|x| = \eta < 1$. To do this, we need two ingredients. 
First, the trace-free Ricci tensor $T$ transforms fairly simply under
conformal changes: if $g = \phi^{-2} g_0$ are any two conformally
related metrics, with trace-free Ricci tensors $T$ and $T_0$, then
\[
T = T_0 + (n-2)\phi^{-1}\left( Dd\phi - \frac{\Delta_0 \phi}{n} g_0\right).
\]
Next, the conformal Killing vector fields on the sphere constitute a
vector space of dimension $(n+1)(n+2)/2$. A spanning set is given by
the following four basic types of conformal Killing fields:
\begin{itemize}
\item $X^{(b)} = \sum b_i \del_{x_i}$, the generators of parabolic
motions fixing infinity,
\item ${\cal D} = \sum x_i \del_{x_i}$, the generator of dilation,
\item $R^{(b,c)} = \sum \left( (b \cdot x)c_i - (c \cdot x)b_i\right)
\del_{x_i}$, the generators of rotations, and
\item $Y^{(b)} = \sum \left( (b \cdot x) x_i - \ha |x|^2 b_i \right)
\del_{x_i}$, the generators of parabolic motions fixing zero.
\end{itemize}

We shall compute relative to the Euclidean metric $\delta$.
The transformation rule for the trace-free Ricci tensor yields
\begin{equation}
T = \frac{2n}{n-2}u^{-2} du\cdot du - 2 u^{-1}Ddu - \left(
\frac{2}{n-2}u^{-2}|\nabla u|^2 + \frac{n-2}{2}u^{\frac{4}{n-2}}
\right)\delta.
\label{eq:tfreuc}
\end{equation}
Next, the hypersurfaces $\Sigma_{\eta}$ have unit normal and 
volume form
\[
\nu =  u^{\frac{-2}{n-2}} \frac{x}{|x|} = u^{\frac{-2}{n-2}}|x|{\cal D},
\qquad \mbox{ and } \quad d\sigma = u^{\frac{2(n-1)}{n-2}}|x|^{n-1}\,d\theta.
\]
Thus we need to evaluate
\begin{eqnarray}
\frac{n-2}{2(n-1)}\int_{|x| = \eta} 
T(X,\nu)\,d\sigma =  \int_{|x| = \eta} 
|x|^{n-2}\left\{ \frac{n}{n-1} (X \cdot \nabla u)({\cal D} \cdot \nabla u) 
- \frac{n-2}{n-1}u\,(Ddu)(X,{\cal D}) 
\right. \nonumber \\
\left. - \left(\frac{1}{n-1}|\nabla u|^2
+ \frac{(n-2)^2}{4(n-1)}u^{\frac{2n}{n-2}}\right)\left(X \cdot {\cal D}
\right) \right\}\,d\theta.
\label{eq:pgen}
\end{eqnarray} 
Denote by $S$ the symmetric $2$-tensor appearing on the
right here (and including the factor $|x|^{n-2}$). 

The main considerations in the computations below 
involve homogeneity. Let us say that a term is of
order $j$ if it is a sum of products of terms
homogeneous of order $j$ and functions periodic
in $\log |x|$. Thus, from the refined asymptotics 
of $u$ we see that 
\[
S = S_{(-2)} + S_{(-1)} + S',
\]
where $S_{(j)}$, $j = -2, -1$, is a symmetric bilinear form, the
coefficients of which (in a standard Euclidean basis)
are of order $j$. All coefficients in the remainder term $S'$ 
are bounded. Each of the conformal fields $X$ listed above
is homogeneous, and so the corresponding term $S_{(j)}(X,{\cal D})$
is of order $j + \ell + 1$ if the coefficients of $X$ are
homogeneous of order $\ell$. Since the integral itself
is independent of the radius $\eta$, the invariant
in each case, coincides with the integral of the term of order zero. 

The easiest case is when $X = Y^{(b)}$, since it has coefficients 
homogeneous of order $2$, and so $S(X,{\cal D}) = O(|x|)$.
We conclude that
\begin{equation}
{\cal P}(Y^{(b)},g) = 0.
\label{eq:pohpar2}
\end{equation}

The other cases require more work. First we make some 
general calculations. Using the refined asymptotics
of $u$ always, we first see that
\begin{equation}
{\cal D} \cdot \nabla u = |x|^{\frac{2-n}{2}}\left\{
-(v_\e' + \frac{n-2}{2}v_\e) + 
(a \cdot x)\left[ -v_\e' + \frac{n-2}{2}v_\e
- \frac{n(n-2)}{4}v_\e^{\frac{n+2}{n-2}}\right] + \ldots \right\}.
\label{eq:rad}
\end{equation}
Next, by definition of the Hessian of a function, 
\[
Ddu(X,{\cal D}) = X({\cal D}u) - (\nabla_X{\cal D})u.
\]
Since $\nabla_X{\cal D} = X$, this reduces to
$X({\cal D}\cdot \nabla u - u)$. Thus we will need that
\begin{equation}
{\cal D}\cdot \nabla u - u = |x|^{\frac{2-n}{2}}
\left\{ -(v_\e' + \frac{n-2}{2}v_\e) - \frac{n(n-2)}{4}(a \cdot x)
v_\e^{\frac{n+2}{n-2}} + \ldots \right\}.
\label{eq:rad2}
\end{equation}
Finally,
\begin{equation}
u^{\frac{2n}{n-2}} = |x|^{-n}\left(v_\e^{\frac{2n}{n-2}}
+ \frac{2n}{n-2}(a \cdot x)(-v_\e' + \frac{n-2}{2}v_\e) + \ldots \right).
\label{eq:hp}
\end{equation}

First consider $X = {\cal D}$.  Since $({\cal D} \cdot \nabla u)^2
= |x|^2 |\nabla u|^2$, the first and third terms combine to give 
\begin{eqnarray}
({\cal D} \cdot \nabla u)^2 = |x|^{2-n}\left(
(v_\e' + \frac{n-2}{2}v_\e)^2 + (a \cdot x)
\left( (v_\e')^2 - \frac{(n-2)^2}{4} + 
\right. \right. \nonumber \\
\left. \left. \frac{n(n-2)}{2}v_\e^{\frac{n+2}{n-2}}(v_\e' + 
\frac{n-2}{2}v_\e)\right) + \ldots \right).
\label{eq:radrad}
\end{eqnarray}
Combining this with the appropriate multiples of ${\cal D}$ applied to 
(\ref{eq:rad2}) and $|x|^2$ multiplied by (\ref{eq:hp}), we 
find that
\[
S({\cal D},{\cal D}) = H(\e)(1 + n(a \cdot x) + \ldots).
\]
Integrating on $S^{n-1}$, and taking the limit as $\eta \rightarrow 0$,
we get 
\begin{equation}
{\cal P}({\cal D},g) = \omega_{n-1}H(\e).
\label{eq:pohrad}
\end{equation}

Next, consider the case $X = R^{(b,c)}$ (or simply $R$, for short).
The third and fourth terms of the integrand involve the
inner product of $R$ with ${\cal D}$, and hence vanish. 
Also, $R$ annihilates any function of $|x|$. Thus, the
components of order zero in each of the remaining terms 
$(R \cdot \nabla u)({\cal D} \cdot \nabla u)$ and $R ({\cal D} 
\cdot \nabla u - u)$ are of the form $F(|x|)(e \cdot x)$ for some 
vectors $e$. These integrate to zero, and so 
\begin{equation}
{\cal P}(R^{(b,c)},g) = 0.
\label{eq:pohrot}
\end{equation}

For the final case, when $X = X^{(b)}$, we decompose $X$
into the orthogonal sum of two vectors $X^{(1)} + X^{(2)}$,
where 
\[
X^{(1)} = \frac{b \cdot x}{|x|^2} {\cal D}
\]
is the radial component. Then
\[
S(X^{(1)},{\cal D}) = \frac{b \cdot x}{|x|^2}S({\cal D},{\cal D}).
\]
Using the expansion above, this becomes
\[
\frac{b \cdot x}{|x|^2}H(\e) + n\frac{(a \cdot x)(b \cdot x)}{|x|^2}
H(\e) + \ldots.
\]
The first term here integrates to zero, while for the second we
use the identity 
\[
\int_{|x| = 1} (a \cdot x) (b \cdot x)\, d\theta = \frac{1}{n}\omega_{n-1}
(a \cdot b),
\]
to conclude that its integral equals
\[
\omega_{n-1}(a \cdot b)H(\e).
\]

For the final term, $S(X^{(2)},{\cal D})$, we continue with
the same methods, noting that since $X^{(2)}$ is orthogonal
to the radial direction, the third and fourth terms once
again vanish. A straightforward calculation, using the
same sorts of parity considerations, shows that the integral 
of the remaining two terms is exactly the same
as for the first component. Putting these together, we
get, at last, that
\begin{equation}
{\cal P}(X^{(b)},g) = \frac{2n}{n-2}\omega_{n-1}H(\e) (a \cdot b).
\label{eq:pohpar1}
\end{equation}

Notice that one consequence of these calculations is that
for any solution $u$ of (\ref{eq:1.1}) with an isolated
singularity at the origin, we can recover the values
of the parameters $\e$ and $a$ in its refined expansion
from the Poho\^{z}aev invariants. 

\subsection{The balancing formul\ae}
We now put the information from the last subsection together as 
follows. Let $u$ correspond to any element of $\Mk$. We regard $u$ as 
a function on $\RR^n$ rather than on the sphere $S^n$ for simplicity.
Assume that $u$ is singular at the points $\{p_1, \ldots, p_k\}$,
and has Fowler and translation parameters at $p_j$ given
by $\e_j$, $T_j$ and $a_j$, respectively. Let $X$ be any one of the
conformal Killing vector fields listed earlier. We compute the
invariants for the homologically trivial hypersurfaces
\[
\Sigma_\eta = \cup_{j=1}^k B(p_j,\eta) \cup B(0,1/\eta)
\]
by letting $\eta$ tend to zero. To do this, first note that 
$u = O(|x|^{2-n})$ as $|x|$ tends to infinity (since $\infty$ 
corresponds to a regular point of $u$ on the sphere), and
thus the integral around the sphere of radius $1/\eta$ tends
to zero with $\eta$. As for the integrals around the other
components of $\Sigma_\eta$, we decompose $X$ at each
$p_j$ as a sum of these four types of basis vector fields translated 
by $p_j$. Thus
\begin{itemize}
\item $X^{(b)} = \sum b_i \del_{x_i}$
\item $D = \sum \big(x_i - p_{j,i}\big)\del_{x_i} + \sum p_{j,i}\del_{x_i}$
\item $R^{(b,c)} = \sum \big( (b \cdot p_j) c_i - (c \cdot p_j) b_i \big)
\del_{x_i} + \sum \big( (b \cdot (x - p_j))c_i - (c \cdot (x - p_j)b_i
\big) \del_{x_i}$
\item $Y^{(b)} = \sum \big( (b \cdot p_j) p_{ji} - \ha |p_j|^2 b_i \big)
+ \sum \big( (b \cdot (x - p_j))p_{ji} - (p_j \cdot (x - p_j))b_i \big)
\del_{x_i} + \sum (b \cdot p_j)(x_i - p_{ji})\del_{x_i} + 
O(|x-p_j|^2).$
\end{itemize}
Here $p_{ji}$ is the $i^{\mbox{\rm th}}$ coordinate of the point $p_j$,
and the last term in the last item is some vector field with
coefficients vanishing quadratically at $p_j$. 

Now apply the computations of the last subsection to evaluate
the limit as $\eta$ tends to zero of the integral around
the sphere $|x-p_j| = \eta$. For $X = X^{(b)}$ we obtain 
$\sum H(\e_j)(b \cdot a_j) = 0$ and since this holds for any vector 
$b \in \RR^n$, we conclude that 
\begin{equation}
\sum_{j=1}^k H(\e_j)a_j = 0.
\label{eq:5.dil}
\end{equation}
Next, with $X = {\cal D}$ we get
\begin{equation}
\sum_{j=1}^k H(\e_j)( (p_j \cdot a_j) + \ha) = 0.
\label{eq:5.rad}
\end{equation}
Next, with $X = R^{(b,c)}$ we get
\begin{equation}
\sum_{j=1}^k H(\e_j) \left( (b\cdot p_j)(c \cdot a_j) - 
(c \cdot p_j)(b \cdot a_j)\right) = 0
\label{eq:5.rot}
\end{equation}
for any vectors $b,c \in \RR^n$. Finally, with $X = Y^{(b)}$ we get
\begin{equation}
\sum_{j=1}^k H(\e_j)\left( (b \cdot p_j)( a_j \cdot p_j + \ha)
- \ha |p_j|^2 (b \cdot a_j) \right) = 0
\label{eq:5.par}
\end{equation}
for any $b \in \RR^n$. These are the global balancing formul\ae\ 
which hold for the parameters associated to any solution
$u \in \Mk$. 

This set of analytic equations does not determine all the parameters 
of the solution $u$, and in particular, gives absolutely no 
information about the `translation parameter' $T_j$ appearing in the 
asymptotics formula at each $p_j$. Still, particularly when 
$k$ is small relative to $n$, these equations do shed more 
light on the global nature of the moduli space. For example, when $k=3$
this set of equations can be solved and leads to an explicit formula 
for the parameters $a_j$ in terms of the $\e_j$ and $p_j$. This
is not possible in general, but nevertheless we still 
have the following result:
\begin{proposition}
There exists a constant $C>0$, depending only on $\Lambda= 
\{p_1, \ldots, p_k\}$, such that for any $u \in {\cal M}_{\Lambda}$, 
the corresponding translation parameters and necksizes satisfy
$|H(\e_j)a_j| \leq C$, $j = 1, \ldots, k$. 
\end{proposition}
We have already seen that ${\cal P}(X^{(b)},g,\del B(p_j,\eta))$
determines $\omega_{n-1}(a_j \cdot b)H(\e_j)$. On the other hand, 
using the universal upper bound of Theorem~\ref{th:3.1} and
elliptic estimates, we also get an {\it a priori} bound for  
this invariant. Since this is true for all $b \in \RR^n$, the
result follows. 

\section{Nondegeneracy of the moduli spaces near their ends}

In this final section we give another application of the asymptotics 
theorem and address the issue of the nondegeneracy of the unmarked 
moduli space $\Mk$. Unlike in the previous sections, we use only the simpler
asymptotics result, not the more refined one. 
We have included this here because the arguments are soft, and not
too different in spirit from some of the ones used above.
This result is an adaptation of one in \cite{MP2}, and our 
desire here is to show its validity beyond the more limited 
setting of that paper.

It is unknown whether degenerate solutions ever exist. If they
do they are quite unstable: it is shown in \cite{MPU1} that under
an arbitrarily small generic change in the conformal class
$[g_0]$, the moduli space becomes smooth. On the other hand,
it is also nontrivial to show that a given solution is 
nondegenerate. This should be easier with explicitly constructed
solutions, but unfortunately this was still impossible to do
with the first-known solutions from \cite{S1}. One construction
of nondegenerate solutions, in a somewhat limited setting,
was given in \cite{MPU2}, and another much more general
one was given in \cite{MP2}. In this last paper it was shown
that given any singular set $\La$, there is a nondegenerate
solution singular at the points of $\La$; these solutions have
very small necksizes (i.e. Fowler parameters). (It was also
shown that for generic configurations $\La$, these solutions
are also nondegenerate in the marked moduli space $\ML$.)

We show here that an argument from \cite{MP2} may be adapted to
prove something slightly weaker, although probably optimal.
Before we state it, we discuss briefly the compactification
theory of these moduli spaces. It was shown in \cite{Po} that 
if $u_\ell$ is any sequence of elements in $\Mk$ such that
the singular points $p_j^{(\ell)}$ are bounded away from one another
and all Fowler parameters $\e_j^{(\ell)}$ are bounded away
from zero, then  some subsequence of $u_\ell$ converges to an element $u_\infty \in \Mk$.
This result limits the ways in which noncompactness in $\Mk$
can occur. This was discussed further in \cite{MPU1} and \cite{MPo}, 
where it was shown that if $u_\ell$ is a sequence in $\Mk$ with
one or both of these restrictions not satisfied, then it is
possible find conformal transformations $F_\ell$ such that
some subsequence of $F_\ell^* u_\ell$ converges to an element $u_\infty'$
in some ${\cal M}_{k'}$ with $k' < k$, or else converges to zero
uniformly on compact sets. (To make things consistent here,
we let ${\cal M}_0$ denote the set of pullbacks of the standard
(smooth) metric on $S^n$ by conformal transformations.)
This result states, then, that $\Mk$ may be compactified
by adding to its ends certain subsets of moduli spaces
${\cal M}_{k'}$ of solutions with fewer or no singular points. 

Finally then we can state our result.
\begin{proposition} Let $u_\ell \in \Mk$ be any sequence of elements
such that the singular points $p_j^{(\ell)}$ stay bounded
away from one another. Suppose that this sequence converges to
$u_\infty \in {\cal M}_{k'}$ with $k' \le k$. Then either $u_\infty$ 
is a degenerate solution in ${\cal M}_{k'}$ or else $u_\ell$ is 
nondegenerate for sufficiently large $\ell$. 
\end{proposition}
For simplicity in the notations, we will assume that the singularities
$\{ p^{(\ell)}_1, \ldots, p^{(\ell)}_k \}$ do not depend on $l$.

The proof is by contradiction. First let us apply the refined
asymptotics theorem as follows. Choose small balls $B(p_j,\rho)$ 
which are disjoint from one another and such that we may write each 
$u_\ell$ as a sum of two functions 
\[
u_\ell = u_{\ell}^{\be,\ba, \bar{T}} + w_\ell.
\]
Here $\be = \{\e_1^{(\ell)}, \ldots, \e_k^{(\ell)}\}$,
$\bar{T} = \{T_1^{(\ell)}, \ldots, T_k^{(\ell)}\}$ and
 $\ba = \{a_1^{(\ell)}, \ldots, a_k^{(\ell)}\}$
are the Fowler and translation parameters at each $p_j$ for $u_\ell$
and $u_\ell^{\be,\ba \bar{T}}$ is a function agreeing with $u_\ell$ outside
the balls $B(p_j,\rho)$ and equalling the model deformed
Fowler solution (relative to the background spherical metric)
$u_{\e_j,a_j, T_j}$ in $B(p_j,\ha \rho)$; finally, $w_\ell \in
{\cal C}^{2,\al}_{\tilde{\gamma}}$ for some fixed $\tilde{\gamma} > (4-n)/2$ close to $(4-n)/2$. 
Assume that at least some of the $\e_j^{(\ell)}$ tend to zero;
otherwise the theorem is trivial. Relabel the points so that,
after passing to a subsequence, $\e_j^{(\ell)}$ tends to zero
for $k'+1 \le j \le  k$ while $\e_j^{(\ell)}$ converges to
some nonzero values $\e_j$ for $j \le k'$. Finally, assume
that there exists, for each $\ell$, a function
$\phi_\ell \in {\cal C}^{2,\al}_{\tilde{\gamma}}$ such that
${\Bbb L}_\ell \phi_\ell = 0$, where ${\Bbb L}_\ell$ is the Jacobi
operator at $u_\ell$. 

Although it is not literally true, we shall assume that the singular
points $p_j$ do not vary with $\ell$. Since we are assuming that they
stay a bounded distance away from one another, we could transform
to this case where $\Lambda$ is fixed by a convergent set of
diffeomorphisms of the sphere, but the only effect this would
have would be to complicate notation.

Normalize $\phi_\ell$, multiplying it by a suitable constant, so that
$\sup d(y)^{-\tilde{\gamma}}|\phi_\ell(y)| = 1$, where $d(y)$ is the distance of the
point $y$ from the singular set $\Lambda$ in the spherical metric.  
Choose a point $y_\ell \in \SL$ realizing this supremum, i.e. such
that $d(y_\ell)^{-\tilde{\gamma}}|\phi_\ell(y_\ell)| = 1$.  

If some subsequence of the $y_\ell$ converges to a point
$y_0 \in \SL$, then we may extract a subsequence of the $v_\ell$ 
converging to an element $v' \in S^n \setminus \Lambda'$, where $\Lambda' =
\{p_{1}, \ldots, p_{k'}\}$ and also so that $\phi_\ell$ converges 
to a nontrivial function $\phi$ on $S^n \setminus \Lambda'$.
Clearly, $|\phi| \le d(y)^{\tilde{\gamma}}$, and also ${\Bbb L}'\phi = 0$ weakly on all
of $S^n$, where ${\Bbb L}'$ is the Jacobi operator at $v'$. 
Since $\tilde{\gamma} > (4-N)/2$ and $v'$ is smooth at the points $p_j$, $j > k'$, 
it follows from a standard removable singularities theorem that $\phi$ is 
actually smooth across these points. Thus $\phi \in {\cal C}^{2,\al}_{\tilde{\gamma}}(S^n \setminus \Lambda')$, and so we have shown that the limiting
solution $v' \in {\cal M}_{k'}$, is degenerate.

If, on the other hand, some subsequence of the $y_\ell$ converges
to one of the points $p_j$, then it is more convenient to
transform the problem, using its conformal equivariance, 
to one on a cylinder before proceeding further. 
First, choose a function $A$ on $S^n \setminus \{p_{j},q\}$ (for any second
point $q$) such that $A^{-\ppp}g_0$ is the product
metric $g_C = \frac{n-2}{n}(dt^2 + d\theta^2)$ on the cylinder 
$C = \RR \times S^{n-1}$, with $t = \infty$
corresponding to $p_j$ and $t = -\infty$ corresponding to $q$. 
On $C$, the function $A$ is simply a multiple of $(\cosh t)^{\frac{n-2}{2}}$,
and on $S^n$ is of the order $\mbox{dist}(y,\{p_j,q\})^{\frac{n-2}{2}}$.  
The solutions $v_\ell$ on $S^n$ correspond to solutions $Av_\ell$ on $C$.
Since the metrics $g_C$ and $g_0$ both have scalar curvature
$n(n-1)$, a straightforward calculation (cf. \cite{MPU2}) shows
that the linearized scalar curvature operators ${\Bbb L}_\ell$ on $S^n$ at $v_\ell$ 
and ${\Bbb L}_{C,\ell}$ on $C$ at $Av_\ell$ satisfy the same conformal
equivariance property as the conformal Laplacians for these two metrics, 
namely 
\begin{equation}
{\Bbb L}_{C,\ell}(A\phi) = A^{\frac{N+2}{N-2}}{\Bbb L}_\ell \phi,
\label{eq:10.4}
\end{equation}
for any function $\phi$. (This depends strongly on the fact that
both metrics have the same scalar curvature.)

Using this transformation, we now have a sequence $Av_\ell$ 
of solutions of the nonlinear equation as well as a sequence
$A\phi_\ell$ of solutions of the Jacobi operator at $Av_\ell$
on $C$. For simplicity we relabel these functions and the operator
by $v_\ell$, $\phi_\ell$ and ${\Bbb L}_\ell$ again. Let $y = (t,\theta)$ 
denote the variable on the cylinder and define $\gamma = \tilde{\gamma} + 
(n-2)/2$, so that $\gamma > 1$, close to $1$. Then 
\begin{equation}
\sup d(y)^{-\gamma}|\phi_\ell| = 1,
\label{eq:10.5}
\end{equation}
where $d(y)$ is once again a smoothed distance function to
the singular points in $\Lambda \setminus \{p_j\}$, transplanted
to $C$, in some large compact set and equal to
$\mbox{sech\,} t$ outside this neighbourhood. By (\ref{eq:10.5}),
$\phi_\ell$ decays at both ends of the cylinder. 

As before, let the supremum in (\ref{eq:10.5}) be attained
at the point $y_\ell = (t_\ell,\theta_\ell)$. By assumption,
$t_\ell \rightarrow \infty$. Translating back by $t_\ell$ and
renormalizing the solution, we find yet another sequence
of solutions, which we again call $w_\ell$, attaining their maximum
at $t = 0$, and which solve the translated equation, which we
again write as ${\Bbb L}_\ell \phi_\ell = 0$. Here ${\Bbb L}_\ell$
is the Jacobi operator at $v_\ell(t+t_\ell,\theta)$.
As usual, some subsequence of the $\phi_\ell$ converge to a nontrivial
solution $\phi$ of the limiting equation $L \phi = 0$, and
$\phi$ is bounded by $e^{-\gamma t}$ for all $t$. 

There are two cases to consider. In the first, $p_j$ is one of 
the singular points for which $\e_j^{(\ell)}$ tends to zero.
Here there are two subcases, depending on whether $v_\ell(t_\ell,
\theta_\ell)$ tends to zero or not. If it does tend to zero, then 
$\phi$ satisfies the equation 
\[
\frac{n}{n-2}\left(\del_t^2 + \Delta_{\th}\right)\phi - 
\frac{n(n-2)}{4}\phi = \frac{n}{n-2}\left( \del_t^2 + \Delta_{\th} - 
\frac{(n-2)^2}{4}\right) \phi = 0.
\]
Decomposing $\phi$ into its $\Lap_\th$ eigencomponents, we see that 
any eigencomponent $\phi_j$ is a sum of exponentials,
$\phi_j = a_j^+ e^{\mu_j t} + a_j^- e^{- \mu_j t}$. Since $\phi$
decays as $t \rightarrow +\infty$, $a_j^+ = 0$. But then it is
clear that no function of the form $e^{-\mu_j t}$ can be bounded
for all $t$ by $e^{-\gamma t}$ unless $\gamma = \mu_j$, which is not
the case, so we arrive at a contradiction.
In the other subcase, $v_\ell(t_\ell,\theta_\ell)$ does not tend to zero. 
Translating by a fixed finite amount, we may assume that $v_\ell$ 
tends to the function $(\cosh t)^{(2-N)/2}$, and hence, after pulling out the
superfluous constants, that the limiting function $\phi$ satisfies
\[
\left(\del_t^2 + \Delta_{\th} - \frac{(n-2)^2}{4} + \frac{n^2 - 4}{4}
\mbox{sech\,}^2 t \right) \phi = 0.
\]
Again separate $\phi$ into its eigencomponents $\phi_j$. Then
\[
\del_t^2 \phi_j - \left(\frac{(n-2)^2}{4} + \lambda_j\right) \phi_j +
\frac{n^2 - 4}{4} \mbox{sech\,}^2 t\, \phi_j = 0.
\]
For $j = 0$ the indicial roots of this equation at both $\pm \infty$
are $\pm (n-2)/2$, for $j = 1$ they are $\pm n/2$, and for 
$j > 1$ the indicial roots are all $ \ge (n+2)/2$. 

The components $\phi_j$ with $j > 1$ are easy to eliminate. In 
fact, these $\phi_j$ must decay faster than $e^{\pm (n+2)|t|/2}$
at $\pm \infty$, so we may multiply the equation satisfied by $\phi_j$ 
and integrate by parts to obtain
\[
\int_{-\infty}^{\infty} (\del_t w_j)^2 + \left(\lambda_j
+ \frac{(n-2)^2}{4}\right) \phi_j^2  - \frac{n^2 - 4}{4}\mbox{sech\,}^2 
t\, \phi_j^2 \,dt =0.
\]
Since $\lambda_j \ge 2n$, the integrand is nonnegative, hence $\phi_j = 0$. 

For the remaining cases, when $j = 0,1$, the indicial roots
at $\pm \infty$ are less than $(n+2)/2$ in absolute
value. On the other hand, to check unmarked nondegeneracy it
suffices to use any $ \gamma_{\e_j,n+1} > \gamma > 1$.
But because $\e_j$ tends to zero, this upper limit tends to $(n+2)/2$. 
Thus if we choose $\gamma$ in the range $(n/2, (n+2)/2)$ then 
both $\phi_{0}$ and $\phi_1$ decay less quickly than $e^{-\gamma t}$ as 
$t \rightarrow +\infty$, which implies that these $\phi_j$ too must vanish.
This is a contradiction. 

The final case to consider is when $\e_j^{(\ell)}$ does not
converge to zero, so that $v_\ell$ is converging to some
Fowler solution $v_\e$. In this case $\phi_\ell$ converges to 
a solution of ${\cal L}_\e \phi = 0$ which satisfies $|\phi| \le
Ce^{-\gamma t}$. But we have shown that all solutions of this
equation are sums of terms for each eigencomponent of $\Lap_\th$ 
which satisfy bounds $|\psi_{\e,j}^\pm| e^{\mp \gamma_{\e,j} t}
\le C$. Clearly this is incompatible with the previous bound,
so we arrive here too at a contradiction.

We have shown, finally, that unless it is converging
to a degenerate solution in some ${\cal M}_{k'}$, $v_\ell$
is nondegenerate in $\Mk$. 

It may well seem disappointing that we can not exclude degeneracy
near any end, but as remarked earlier, this form of the result
is probably optimal. However, in certain cases we can deduce
nondegeneracy without restriction; this is due to the fact that only
the first of the several cases treated in the proof did not necessarily
lead to a contradiction.  However, for example, if the solutions 
$v_\ell$ are converging to zero uniformly on compact sets of $\SL$, 
then they must be nondegenerate for sufficiently large $\ell$. This 
is because the equation satisfied by the limiting Jacobi field $\phi$
is $(\Lap_{S^n} - (n-2)^2/4)\phi = 0$ which has no nontrivial
solution. This is the case studied in \cite{MP2}. Another case where
we can deduce nondegeneracy for sufficiently large $\ell$ is when
all but two of the singular points $p_j$ disappear in the limit,
i.e. when $k' = 2$. This is because any Fowler solution is
nondegenerate. On the other hand, it is not possible, using this
argument, to deduce nondegeracy for the solutions constructed 
in \cite{S1}, for there $v_\ell$ converges to the constant function
$1$ on $S^n$, and the sphere is degenerate.

\end{document}